\documentclass[number,citesort,seceqn,dvips]{arxbj}
\usepackage{upgreek,multirow}
\usepackage{graphicx}

\aid{0}
\volume{18}
\issue{1}
\pubyear{2012}
\firstpage{177}
\lastpage{205}
\doi{10.3150/10-BEJ331}

\makeatletter
\newcommand{\bzero}{{\mathbf0}}
\newcommand{\bff}{\mathbf{f}}
\newcommand{\bI}{\mathbf{I}}
\newcommand{\bM}{\mathbf{M}}
\renewcommand{\hm}{{\hat{m}}}
\newcommand{\hg}{{\hat{g}}}
\newcommand{\hPi}{{\hat{\Pi}}}
\newcommand{\hQ}{{\hat{Q}}}
\newcommand{\tm}{{\tilde{m}}}
\newcommand{\tg}{{\tilde{g}}}
\newcommand{\tbm}{{\tilde\mathbf{m}}}
\newcommand{\hbr}{{\hat\mathbf{r}}}
\newcommand{\bZ}{\mathbf{Z}}
\newcommand{\bN}{\mathbf{N}}
\newcommand{\hbM}{{\hat\mathbf{M}}}
\newcommand{\bV}{\mathbf{V}}
\newcommand{\bl}{\mathbf{l}}
\newcommand{\bh}{\mathbf{h}}
\newcommand{\bg}{\mathbf{g}}
\newcommand{\bX}{\mathbf{X}}
\newcommand{\bu}{\mathbf{u}}
\newcommand{\bw}{\mathbf{w}}
\newcommand{\betam}{\bolds\beta}
\newcommand{\deltam}{\bolds\delta}
\newcommand{\bmu}{\bolds\mu}
\newcommand{\etam}{\bolds\eta}
\newcommand{\gamm}{\bolds\gamma}
\newcommand{\Psim}{\bolds\Psi}
\newcommand{\hPsim}{{\hat{\Psim}}}

\newtheorem{theorem}{Theorem}
\newproclaim{BA}{Backfitting algorithm}

\def\cH{{\mathcal H}}
\newcommand{\cal}{\mathcal}
\newcommand{\overset}{\stackrel}
\newproclaim{assumptions}{Assumptions}
\makeatother

\begin{document}
\begin{frontmatter}

\title{Projection-type estimation for varying
coefficient regression models}
\runtitle{Varying coefficient models}

\begin{aug}
\author[a]{\fnms{Young K.} \snm{Lee}\thanksref{a}\ead[label=e1]{youngklee@kangwon.ac.kr}},
\author[b]{\fnms{Enno} \snm{Mammen}\thanksref{b}\ead[label=e2]{emammen@rumms.uni-mannheim.de}} \and
\author[c]{\fnms{Byeong U.} \snm{Park}\corref{}\thanksref{c}\ead[label=e3]{bupark@stats.snu.ac.kr}}
\runauthor{Y.K. Lee, E. Mammen and B.U. Park}
\address[a]{Department of Statistics, Kangwon National University,   Chuncheon  200-701,  Korea.\\
\printead{e1}}
\address[b]{Department of Economics, University of Mannheim, L7, 3-5, 688131 Mannheim, Germany.\\
\printead{e2}}
\address[c]{Department of Statistics, Seoul National University, Seoul 151-747, Korea.\\
\printead{e3}}
\end{aug}

\received{\smonth{5} \syear{2010}}
\revised{\smonth{9} \syear{2010}}

%
\begin{abstract}
In this paper we introduce new estimators of the coefficient functions
in the varying coefficient regression model. The proposed estimators
are obtained by projecting the vector of the full-dimensional
kernel-weighted local polynomial estimators of the coefficient
functions onto a Hilbert space with a suitable norm. We provide a
backfitting algorithm to compute the estimators. We show that the
algorithm converges at a geometric rate under weak conditions. We
derive the asymptotic distributions of the estimators and show that the
estimators have the oracle properties. This is done for the general
order of local polynomial fitting and for the estimation of the
derivatives of the coefficient functions, as well as the coefficient
functions themselves. The estimators turn out to have several
theoretical and numerical advantages over the marginal integration
estimators studied by Yang, Park, Xue and H\"{a}rdle [\textit
{J. Amer. Statist. Assoc.} \textbf{101} (2006) 1212--1227].
\end{abstract}

%
\begin{keyword}
\kwd{kernel smoothing}
\kwd{local polynomial regression}
\kwd{ marginal integration}
\kwd{oracle properties}
\kwd{smooth backfitting}
\kwd{varying coefficient models}
\end{keyword}

\end{frontmatter}

\section{Introduction}\label{sec1}

In this paper we consider a varying coefficient regression model
proposed by Hastie and Tibshirani \cite{Hastie1993} and studied by
Yang, Park, Xue
and H\"{a}rdle \cite{Yang2006}. The model takes the form $Y^i = m(\bX
^i, \bZ
^i) + \sigma(\bX^i, \bZ^i) \varepsilon^i$, $i=1, \ldots, n$, where
%
%
\begin{equation}\label{model}
m(\bX, \bZ) = \sum_{j=1}^d m_j(X_j)Z_j,
\end{equation}
$m_j$ are unknown coefficient functions, $\bX^i=(X_1^i, \ldots,
X_d^i)^\top$ and $\bZ^i=(Z_1^i, \ldots, Z_d^i)^\top$ are observed
vectors of covariates, and $\varepsilon^i$ are the error variables
such that \mbox{$E(\varepsilon^i |\bX^i, \bZ^i)=0$} and $\operatorname
{var}(\varepsilon^i |\bX^i, \bZ^i)=1$. We assume that $(\bX^i, \bZ
^i, Y^i)$ for $1 \le i \le n$ are independent and identically
distributed. The model is simple in structure and easily interpreted,
yet flexible, since the dependence of the response variable on the
covariates is modeled in a~nonparametric way. The model is different
from the functional coefficient model of Chen and Tsay \cite
{Chen1993}, Fan and
Zhang \cite{Fan1999}, Cai, Fan and Li \cite{CaiLi2000} and Cai, Fan
and Yao \cite{CaiYao2000}, where
$m_j$ are functions of a single variable, that is, $m(X^i, \bZ^i) =
\sum_{j=1}^d m_j(X^i)Z_j^i$. Fitting the latter model is much simpler
than the model (\ref{model}) since it involves only a~univariate
smoothing across the single variable $X$.

To fit the model (\ref{model}), we may apply the idea of local
polynomial smoothing. To illustrate the difficulty in fitting the
model, suppose that we employ local constant fitting so that we minimize
\[
\sum_{i=1}^n \Biggl[Y^i - \sum_{j=1}^d \theta_j Z_j^i \Biggr]^2
K_h(x_1, X_1^i) \cdots K_h(x_d, X_d^i)
\]
with respect to $\theta_j$, $1 \le j \le d$, to get estimators of
$m_j(x_j)$, $1 \le j \le d$, where $K_h$ is a kernel function. For each
coefficient $m_j$, this yields an estimator which is a function of not
only $x_j$ but also other variables $x_k$, $k \neq j$. The marginal
integration method, proposed and studied by Yang \textit{et al.} \cite
{Yang2006}, is
simply to take the average of $\hat{\theta}_j(X_1^i, \ldots,
X_{j-1}^i, x_j, X_{j+1}^i, \ldots, X_d^i)$ in order to eliminate the
dependence on the other variables.

In this paper we propose a new method for fitting the model (\ref
{model}). The proposed method is to project the vector of the
full-dimensional kernel-weighted local polynomial estimators ($\hat
{\theta}_j$, $1 \le j \le d$, in the above, in the case of local
constant fitting) onto a~space of vectors of functions $f_j\dvtx  \mathbb
{R} \rightarrow\mathbb{R}$, $1 \le j \le d$, with a suitable norm.
Projection-type estimation has been studied in other structured
nonparametric regression models. For example, the smooth backfitting
method was proposed by Mammen, Linton and Nielsen~\cite{Mammen1999} to fit
additive regression models. The same idea was applied to
quasi-likelihood additive regression by Yu, Park and Mammen \cite
{Yu2008} and
to additive quantile regression by Lee, Mammen and Park \cite
{Lee2010}. Some
nonparametric time series models have been proposed with unobserved
factors $Z_j$ that do not depend on the individual but on time; see,
for example,
Connor, Linton and Hagmann \cite{Connor2008}, Fengler, H\"{a}rdle and Mammen
\cite{Fengler2007} and Park, Mammen, H\"{a}rdle and Borak \cite
{Park2009}. In these papers
it has been shown that one can also proceed asymptotically in the
models under consideration, as if the factors would have been observed.
We note that the current problem does not fit into the framework of the
above papers but requires a different treatment. In particular, in the
model (\ref{model}), the functions $m_j$ are not additive components
of the regression function, but they are the coefficients of $Z_j$. For
a treatment of our model we have to exclude the case of constant
$Z_j\equiv1$. In the case of constant $Z_j$, model (\ref{model})
reduces to the additive model. The key element in the derivation of the
theory for our model is to embed the vector of the coefficient
functions into an additive space of vectors of univariate functions and
then to endow the space with a norm where the covariates $Z_j$ enter
with kernel weights.

As far as we know, the marginal integration method has been the only
method to fit the model (\ref{model}). It is widely accepted that the
marginal integration method suffers from the curse of dimensionality.
Inspired by Fan, H\"{a}rdle and Mammen \cite{Fan1998} and others, Yang \textit{et
al.} \cite{Yang2006} tried to avoid the dimensionality problem by using two
different types of kernels and bandwidths. To be more specific,
consider estimation of $m_j$ for a particular~$j$. The method then uses
a kernel, say $L$, and bandwidths, say $b_k$, for the directions of
$x_k$ $(k \neq j)$, which are different from a kernel $K$ and a
bandwidth $h_j$ for the direction of $x_j$. By choosing $b_k \ll h_j$
and taking a higher order kernel $L$, we can achieve the univariate
optimal rate of convergence for the resulting estimator of $m_j$. One
of the main difficulties with the marginal integration method is that
there is no formula available for the optimal choice of the secondary
bandwidths $b_k$. Also, the performance of the method depends crucially
on the choice of the secondary bandwidths $b_k$, as observed in our
numerical study; see Section~\ref{sec5}. Furthermore, the method involves
estimation of a full-dimensional regression estimator, which requires
inversion of a full-dimensional $[(\pi+1)d]\times[(\pi+1)d]$
smoothing matrix, where $\pi$ is the order of local polynomial
fitting. This means that the method may break down in practice in high
dimension.

On the contrary, the proposed method may use bandwidths of the same
order for all directions to achieve the univariate optimal rate of
convergence, and we derive formulas for the optimal bandwidths. The
method requires only one- and two-dimensional smoothing and inversion
of a $(\pi+1) \times(\pi+1)$ matrix which is computed by a
single-dimensional local smoothing. Thus, the proposed method does not
suffer from the curse of dimensionality in practice as well as in
theory. We show that the method has the oracle properties, meaning that
the proposed estimator of $m_j$ for each $j$ has the same first-order
asymptotic properties as the oracle (infeasible) estimator of $m_j$
that uses the knowledge of all other coefficient functions $m_k$, $k
\neq j$. We develop the theory for the method with local polynomial
fitting of general order $\pi\ge0$. Thus, the theory gives the
asymptotic distributions of the estimators of $m_j$, as well as their
derivatives $m_j^{(k)}$, $1 \le k \le\pi$.

There have been several works on a related varying coefficient model
where the coefficients are time-varying functions. These include
Hoover, Rice, Wu and Yang \cite{Hoover1998}, Huang, Wu and Zhou \cite
{Huang2002,Huang2004},
Wang, Li and Huang \cite{Wang2008} and Noh and Park \cite{Noh2010}.
The kernel method of
fitting this model is quite different from, and simpler than, the
method of fitting our model (\ref{model}), since the former involves
only a univariate smoothing across time. Recently, Park, Hwang and Park
\cite{Park2010} considered a testing problem for the model~(\ref
{model}) based
on the marginal integration method.

This paper is organized as follows. In the next section, we describe
the proposed method with local constant fitting and then, in Section~\ref{sec3},
we give its theoretical properties. Section~\ref{sec4} is devoted to the
extension of the method and theory to local polynomial fitting of
general order. In Section~\ref{sec5} we present the results of our numerical
study. In Section~\ref{sec6} we apply the proposed method to Boston Housing
Data. Technical details are contained in the \hyperref[appm]{Appendix}.

\section{The method with local constant fitting}\label{sec2}

Although our main focus is to introduce the method with local
polynomial fitting and to develop its general theory, we start with
local constant fitting since the method is better understood in the
latter setting. Let $Y$ be the response variable, and $\bX=(X_1,\ldots
,X_d)^\top$ and $\bZ=(Z_1, \ldots, Z_d)^\top$ be
the covariate vectors of dimension $d$. Let $\{(\bX^i, \bZ^i, Y^i)\}
_{i=1}^n$ be a random sample drawn from $(\bX, \bZ, Y)$. Assume that
the density $p$ of $\bX$ is supported on $[0,1]^d$. To estimate the
coefficient functions $m_j$ in the model (\ref{model}), we consider a~%
`smoothed' squared error loss. Similar ideas were adopted for additive
regression by Mammen \textit{et al.} \cite{Mammen1999} and for quasi-likelihood additive
regression by Yu \textit{et al.} \cite{Yu2008}.

Let $K$ be a nonnegative function, called a \textit{base kernel}. To
define a smoothed squared error loss, we use a boundary corrected
kernel, as in Mammen \textit{et al.} \cite{Mammen1999} and Yu~\textit{et al.}~\cite
{Yu2008}, which is
defined by
\[
K_g(u,v) = \biggl[\int_0^1 K \biggl(\frac{w-v}{g} \biggr)\,\mathrm{d}w
\biggr]^{-1} K \biggl(\frac{u-v}{g} \biggr) I(u,v \in[0,1]).
\]
Suppose that we use different bandwidths for different directions. Let
$\bh= (h_1, \ldots, h_d)$ be the bandwidth vector. For simplicity, we
focus on the case where we use a product kernel of the form $K_\bh(\bu
, \mathbf{v}) = \prod_{j=1}^d K_{h_j}(u_j, v_j)$. We may use a more
general multivariate kernel, but this would require more involved
notation and technical arguments. The proposed estimator of $\mathbf
{m}\equiv(m_1, \ldots, m_d)^\top\dvtx  \mathbb{R}^d \rightarrow\mathbb
{R}^d$, where $m_j(\mathbf{x}) = m_j(x_j)$, is defined to be the
minimizer of
\[
L(\bff) = \int n^{-1} \sum_{i=1}^n\Biggl [Y^i - \sum_{j=1}^d
f_j(x_j) Z_j^i \Biggr]^2 K_\bh(\mathbf{x}, \bX^i)\,\mathrm{d}\mathbf{x}
\]
over $\bff= (f_1, \ldots, f_d)^\top$ with $L(\bff) < \infty$. Here
and hereafter, integration over $\mathbf{x}$ is on $[0,1]^d$. Define
$\hbM
(\mathbf{x}) = n^{-1}\sum_{i=1}^n K_\bh(\mathbf{x}, \bX^i) \bZ^i
\bZ
^{i\top}$. Then, $L(\bff) < \infty$ is equivalent to $\int\bff
(\mathbf{x})^\top\hbM(\mathbf{x}) \bff(\mathbf{x}) \,\mathrm
{d}\mathbf{x}< \infty$. The
function space that arises in the minimization problem is\looseness=1
\[
\cH(\hbM) = \{\bff\in L_2(\hbM)\dvtx  f_j(\mathbf{x})=g_j(x_j) \mbox{ for
a function $g_j\dvtx \mathbb{R} \rightarrow\mathbb{R}$}, 1 \le j \le d\},
\]\looseness=0
where $L_2(\hbM)$ denotes a class of function vectors $\bff$ defined by
\begin{eqnarray*}
L_2(\hbM) &=& \biggl\{\bff\dvtx  \bff(\mathbf{x})=(f_1(\mathbf{x}), \ldots,
f_d(\mathbf{x}))^\top\mbox{ for
some functions }f_j\dvtx \mathbb{R}^d \rightarrow\mathbb{R} \\
&& \hphantom{\biggl\{}\mbox{and }\int\bff(\mathbf{x})^\top\hbM(\mathbf{x}) \bff
(\mathbf{x}) \,\mathrm{d}\mathbf{x}< \infty\biggr\}.
\end{eqnarray*}

The spaces $L_2(\hbM)$ and $\cH(\hbM)$ are Hilbert spaces equipped
with a (semi)norm $\|\cdot\|_\hbM$, defined by
\[
\|\bff\|_\hbM^2 = \int\bff(\mathbf{x})^\top\hbM(\mathbf{x})
\bff(\mathbf{x})\,\mathrm{d}\mathbf{x}.
\]
Let $\bM(\mathbf{x}) = E(\bZ\bZ^\top| \bX=\mathbf{x}) p(\mathbf
{x})$. Since
$\|\bff\|_\hbM^2$ converges to
\[
\|\bff\|_\bM^2 \equiv\int\bff(\mathbf{x})^\top\bM(\mathbf{x})
\bff
(\mathbf{x})\,\mathrm{d}\mathbf{x}
\]
in probability under certain conditions, the corresponding Hilbert
spaces in the\vspace*{1pt} limit are~$L_2(\bM)$ and $\cH(\bM)$, which are defined
as $L_2(\hbM)$ and $\cH(\hbM)$, respectively, with $\hbM$ being
replaced by $\bM$. Here, we note that $\|\cdot\|_\bM$ becomes a norm
if we assume that
%
%
\begin{equation}\label{concurvity}
\bff(\bX)^\top\bZ=0 \mbox{ almost surely} \mbox{ implies }
\bff= \bzero.
\end{equation}
In fact, the assumption (\ref{concurvity}) is known to be a sufficient
condition for avoiding \textit{concurvity}, as termed by Hastie and
Tibshirani \cite{Hastie1990}, an analog of collinearity in linear
models. If the
assumption does not hold, then the $m_j$ are not identifiable. This is
because, for $\bff$ such that $\bff(\bX)^\top\bZ= 0$ almost
surely, we have
\[
E(Y | \bX, \bZ) = \mathbf{m}(\bX)^\top\bZ= [\mathbf{m}(\bX
)+\bff
(\bX) ]^\top\bZ.
\]
The assumption (\ref{concurvity}) is satisfied if we assume that the
smallest eigenvalue of $E(\bZ\bZ^\top|\allowbreak \bX=\mathbf{x})$ is bounded
away from zero on $[0,1]^d$.

For $\bff\in\cH(\hbM)$, we obtain
\begin{eqnarray*}
L(\bff) &=& \int n^{-1} \sum_{i=1}^n [Y^i - \tbm(\mathbf{x})^\top
\bZ^i ]^2 K_\bh(\mathbf{x}, \bX^i)\,\mathrm{d}\mathbf{x}\\
&& {}+ \int[\tbm(\mathbf{x})- \bff(\mathbf{x}) ]^\top
\hbM(\mathbf{x}) [\tbm(\mathbf{x})- \bff(\mathbf{x}) ] \,\mathrm
{d}\mathbf{x},
\end{eqnarray*}
where $\tbm$ is the minimizer of $L(\bff)$ over $\bff\in L_2(\hbM
)$. It is given explicitly as
%
%
\begin{equation}\label{def-tbm}
\tbm(\mathbf{x}) = \hbM(\mathbf{x})^{-1} n^{-1} \sum_{i=1}^n \bZ
^i Y^i
K_\bh(\mathbf{x}, \bX^i).
\end{equation}
Thus, the proposed estimator ${\hat\mathbf{m}}= (\hm_1, \ldots, \hm
_d)^\top$ can be defined equivalently as the projection of $\tbm$
onto $\cH(\hbM)$:
%
%
\begin{equation}\label{def-est}
{\hat\mathbf{m}}= \operatorname{argmin}\limits_{\bff\in\cH(\hbM)} \|
\tbm- \bff\|_\hbM^2.
\end{equation}

By considering the G\^{a}teaux or Fr\'{e}chet derivatives of the
objective function with respect to~$\bff$, the solution ${\hat\mathbf
{m}}$ of the minimization problem (\ref{def-est}) satisfies the
following system of integral equations:
%
%
\begin{equation}\label{backeqn}
0 = \int\hbM_j(\mathbf{x})^\top[\tbm(\mathbf{x}) - {\hat\mathbf
{m}}(\mathbf{x}) ] \,\mathrm{d}\mathbf{x}_{-j}, \qquad  1 \le j \le d,
\end{equation}
where $\hbM_j$ are defined by $\hbM=(\hbM_1, \ldots, \hbM_d)^\top
$ and $\mathbf{x}_{-j} = (x_1, \ldots, x_{j-1}, x_{j+1}, \ldots,
x_d)^\top$. In fact, the system (\ref{backeqn}) turns out to be a
backfitting system of equations. To see this, we define
\begin{eqnarray*}
\tm_j(x_j) &=& {\hat{q}}_j (x_j)^{-1} n^{-1}\sum_{i=1}^n
K_{h_j}(x_j, X_j^i) Z_j^i Y^i,\\
{\hat{q}}_j(x_j) &=& n^{-1}\sum_{i=1}^n K_{h_j}(x_j, X_j^i)
(Z_j^i)^2,\\
{\hat{q}}_{jk}(x_j, x_k) &=& n^{-1}\sum_{i=1}^n K_{h_j}(x_j, X_j^i)
K_{h_k}(x_k, X_k^i)Z_j^i Z_k^i, \qquad  k \neq j.
\end{eqnarray*}
We note that, by definition, $\tm_j\dvtx  \mathbb{R} \rightarrow\mathbb
{R}$ does not equal the $j$th component of $\tbm$, which maps $\mathbb
{R}^d$ to $\mathbb{R}$. We can then see that
\begin{eqnarray*}
\int\hbM_j(\mathbf{x})^\top\tbm(\mathbf{x})\,\mathrm{d}\mathbf
{x}_{-j} &=& \tm
_j(x_j){\hat{q}}_j(x_j),\\
\int\hbM_j(\mathbf{x})^\top{\hat\mathbf{m}}(\mathbf{x})\,\mathrm
{d}\mathbf{x}_{-j} &=&
\hm_j(x_j){\hat{q}}_j(x_j) + \sum_{k=1, \neq j}^d \int\hm
_k(x_k){\hat{q}}_{jk}(x_j,x_k) \,\mathrm{d}x_k.
\end{eqnarray*}
The second formula is obtained by using the following property of the
boundary corrected kernel: $\int K_{h_j}(u_j,v_j) \,\mathrm{d}u_j=1$.
Thus, the
system of equations (\ref{backeqn}) is equivalent to
%
%
\begin{equation}\label{backeqn1}
\hm_j(x_j) = \tm_j(x_j) - \sum_{k=1, \neq j}^d \int\hm_k(x_k)\frac
{{\hat{q}}_{jk}(x_j,x_k)}{{\hat{q}}_j(x_j)} \,\mathrm{d}x_k, \qquad  1 \le j
\le d.
\end{equation}

We emphasize that the proposed method does not require computation of
the full-dimensional estimator $\tbm(\mathbf{x})$ at (\ref
{def-tbm}). It
only requires one- and two-dimensional smoothing to compute $\tm_j$,
${\hat{q}}_j$ and ${\hat{q}}_{jk}$, and involves inversion of ${\hat
{q}}_j$ only. In contrast, the marginal integration method studied by
Yang \textit{et al.} \cite{Yang2006} involves the computation of~$\tbm(\mathbf
{x})$, which
requires inversion of the full-dimensional smoothing matrix $\hbM$.
Thus, in practice, the marginal integration method may break down in
high dimensions where $d$ is large.\looseness=1

We express the updating equations (\ref{backeqn1}) in terms of
projections onto suitable function spaces. This representation is
particularly useful in our theoretical development. We consider $\cH
_j(\hbM)$, $1 \le j \le d$, subspaces of $\cH(\hbM)$ defined by
\[
\cH_j(\hbM) = \{\bff\in L_2(\hbM)\dvtx  f_j(\mathbf{x})=g_j(x_j) \mbox
{ for
a function }g_j\dvtx \mathbb{R} \rightarrow\mathbb{R}, f_k \equiv0
\mbox{ for }k \neq j \}.
\]
With this definition, we have
\[
\cH(\hbM) = \cH_1(\hbM)+ \cdots+ \cH_d(\hbM).
\]
Also, denoting the projection operator onto a closed subspace $\cH$ of
$\cH(\hbM)$ by $\Pi(\cdot| \cH)$ and its $j$th element by $\Pi
(\cdot| \cH)_j$, we get, for $\bff\in L_2(\hbM)$,
%
%
\begin{eqnarray}
\label{form1}
\Pi(\bff| \cH_j(\hbM) )_j &=& {\hat
{q}}_j(x_j)^{-1}\int\hbM_j(\mathbf{x})^\top\bff(\mathbf{x})\,
\mathrm{d}\mathbf{x}_{-j}, \nonumber
\\[-8pt]
\\[-8pt]
\Pi(\bff| \cH_j(\hbM) )_k &=& 0, \qquad  k \neq j.
\nonumber
\end{eqnarray}
In particular, for $\bff\in\cH(\hbM)$, we have
%
%
\begin{equation}\label{form2}
\Pi(\bff| \cH_j(\hbM) )_j = f_j(x_j) + \sum_{k=1,
\neq j}^d \int f_k(x_k)\frac{{\hat{q}}_{jk}(x_j,x_k)}{{\hat
{q}}_j(x_j)} \,\mathrm{d}x_k.
\end{equation}
Furthermore, for $\bff\in\cH_k(\hbM)$,
%
%
\begin{equation}\label{proj}
\Pi(\bff| \cH_j(\hbM) )_j = \int f_k(x_k)\frac
{{\hat{q}}_{jk}(x_j,x_k)}{{\hat{q}}_j(x_j)} \,\mathrm{d}x_k, \qquad  j \neq k.
\end{equation}

For ${\hat\mathbf{m}}\in\cH(\hbM)$, let ${\hat\mathbf
{m}}_j(\mathbf{x}) =
(0,\ldots, 0, \hm_j(x_j), 0, \ldots, 0)^\top$ denote the vector
whose $j$th entry equals $\hm_j(x_j)$, the rest being zero. We can
then decompose ${\hat\mathbf{m}}$ as ${\hat\mathbf{m}}= {\hat
\mathbf{m}}_1 +
\cdots+ {\hat\mathbf{m}}_d$. From (\ref{backeqn1}) and (\ref{proj}),
we obtain
%
%
\begin{equation}\label{backeqn2}
{\hat\mathbf{m}}_j = \Pi\Biggl(\tbm- \sum_{k=1, \neq j}^d {\hat\mathbf
{m}}_k \Big| \cH_j(\hbM) \Biggr), \qquad  1 \le j \le d.
\end{equation}
The backfitting equations (\ref{backeqn1}), or their equivalent forms
(\ref{backeqn2}), give the following backfitting algorithm.

\begin{BA*}
With a set of initial estimates $\hm_j^{[0]}$,
iterate for $r=1, 2, \ldots$ the following process: for $1 \le j \le d,$
\begin{eqnarray*}
\hm_j^{[r]}(x_j) &=& \tm_j(x_j) - \sum_{k=1}^{j-1} \int\hm
_k^{[r]}(x_k)\frac{{\hat{q}}_{jk}(x_j,x_k)}{{\hat{q}}_j(x_j)} \,
\mathrm{d}x_k \\
&&{} - \sum_{k=j+1}^{d} \int\hm_k^{[r-1]}(x_k)\frac{{\hat
{q}}_{jk}(x_j,x_k)}{{\hat{q}}_j(x_j)} \,\mathrm{d}x_k
\end{eqnarray*}
or, equivalently,
%
%
\begin{eqnarray}\label{backfit}
{\hat\mathbf{m}}_j^{[r]} &=& \Pi\Biggl(\tbm- \sum_{k=1}^{j-1} {\hat
\mathbf{m}}_k^{[r]} - \sum_{k=j+1}^{d}{\hat\mathbf{m}}_k^{[r-1]}\Big |
\cH
_j(\hbM) \Biggr).
\end{eqnarray}
\end{BA*}

\section{Theoretical properties of the local constant method}\label{sec3}

\subsection{Convergence of the backfitting algorithm}\label{sec3.1}

The theoretical development for the backfitting algorithm (\ref
{backfit}) and for the solution of the backfitting equation (\ref
{backeqn2}) does not fit into the framework of an additive regression
function as in Mammen \textit{et al.} \cite{Mammen1999}. Formally, we get their
model by
taking $Z_j \equiv1$ for all $1 \le j \le d$ in (\ref{model}).
However, for identifiability of $m_j$, we need the assumption that
$E(\bZ\bZ^\top| \bX=\mathbf{x})$ is invertible; see the assumption
(A1) below. Trivially, this assumption does not hold for the additive
model with $Z_j \equiv1$. For our model, we directly derive the
theoretical properties of the algorithm and the estimators by borrowing
some relevant theory on projection operators.

Let $\hPi_j$ denote the projection operator $\Pi(\cdot|  \cH_j(\hbM
))$ and $\Pi_j$ the projection operator $\Pi(\cdot| \cH_j(\bM))$.
Define $\hQ_j = I-\hPi_j$ and $Q_j = I- \Pi_j$; these are the
projection operators onto $\cH_j(\hbM))^\perp$ and $\cH_j(\bM
))^\perp$, respectively. From the backfitting algorithm (\ref
{backfit}), it follows that
%
%
\begin{eqnarray}\label{backfit1}
\tbm- \sum_{k=1}^{j} {\hat\mathbf{m}}_k^{[r]} - \sum
_{k=j+1}^{d}{\hat
\mathbf{m}}_k^{[r-1]} &=& \hQ_j \Biggl(\tbm- \sum_{k=1}^{j-1} {\hat
\mathbf{m}}_k^{[r]} - \sum_{k=j+1}^{d}{\hat\mathbf{m}}_k^{[r-1]}
\Biggr)\nonumber
\\[-8pt]
\\[-8pt]
&=& \hQ_j \Biggl(\tbm- \sum_{k=1}^{j-1} {\hat\mathbf{m}}_k^{[r]} -
\sum_{k=j}^{d}{\hat\mathbf{m}}_k^{[r-1]} \Biggr).
\nonumber
\end{eqnarray}
Define $\hQ= \hQ_d \cdots\hQ_1$. Repeated application of (\ref
{backfit1}) for $j=d, d-1, \ldots, 1$ gives
\[
\tbm- {\hat\mathbf{m}}^{[r]} = \hQ\bigl(\tbm- {\hat\mathbf{m}}^{[r-1]}\bigr).
\]
This establishes
that
%
%
\begin{equation}\label{backfit2}
{\hat\mathbf{m}}^{[r]} = \hQ{\hat\mathbf{m}}^{[r-1]} + \hbr= \sum
_{s=0}^{r-1} \hQ^s \hbr+ \hQ^r {\hat\mathbf{m}}^{[0]},
\end{equation}
where $\hbr= (I-\hQ) \tbm$. If we write $\tbm_j(\mathbf{x}) =
(0,\ldots, 0, \tm_j(x_j), 0, \ldots, 0)^\top$, then $\hPi_j \tbm=
\tbm_j$ so that
%
%
\begin{eqnarray}\label{backfit3}
\hbr&=& (\hPi_d + \hQ_d \hPi_{d-1} + \cdots+ \hQ_d \cdots
\hQ_2 \hPi_1 )\tbm\nonumber
\\[-8pt]
\\[-8pt]
&=& \tbm_d + \hQ_d \tbm_{d-1} + \cdots+ \hQ_d \cdots\hQ_2 \tbm_1.
\nonumber
\end{eqnarray}

Convergence of the backfitting algorithm (\ref{backfit}) depends on
the statistical properties of the operator $\hQ$. Consider the event
${\cal E}_n$, where $\hbr, {\hat\mathbf{m}}^{[0]} \in\cH(\bM)$ and
the norm of the operator $\hQ$ is strictly less than one, that is, $\|
\hQ\| <1$. Here and below, for an operator $F\dvtx  \cH(\bM)
\rightarrow\cH(\bM)$,
\[
\|F\| = \sup\{\|F \bff\|_\bM\dvtx  \bff\in\cH(\bM) , \|\bff\|_\bM
\le1\}.
\]
Then, in that event, $\sum_{s=0}^\infty\hQ^s \hbr$ is well defined
in $\cH(\bM)$ and, by (\ref{backfit2}), ${\hat\mathbf{m}}^{[r]}$
converges to $\sum_{s=0}^\infty\hQ^s \hbr$ as $r$ tends to
infinity. The limit is a solution of the backfitting equation (\ref
{backeqn2}) since the latter is equivalent to ${\hat\mathbf{m}}= \hQ
{\hat\mathbf{m}}+ \hbr$. Furthermore, the solution is unique since
repeated application of ${\hat\mathbf{m}}= \hQ{\hat\mathbf{m}}+
\hbr$
leads to ${\hat\mathbf{m}}=\sum_{s=0}^\infty\hQ^s \hbr$.

Below, we collect the assumptions that make the event ${\cal E}_n$
occur with probability tending to one and state a theorem for the
convergence of the backfitting algorithm~(\ref{backfit}).\vadjust{\goodbreak}

\begin{assumptions*}
\begin{enumerate}[(A11)]
\item[(A1)] $E(\bZ\bZ^\top| \bX=\mathbf{x})$ is continuous and its
smallest eigenvalue is bounded away from zero on $[0,1]^d$.

\item[(A2)] $\sup_{\mathbf{x}\in[0,1]^d} E(Z_j^4 | \bX=\mathbf
{x}) <
\infty$ for all $1 \le j \le d$.

\item[(A3)] The joint density $p$ of $\bX$ is bounded away from zero
and is continuous on $[0,1]^d$.

\item[(A4)] $E|Y|^\alpha< \infty$ for some $\alpha>5/2$.

\item[(A5)] $K$ is a bounded and symmetric probability density
function supported on $[-1,1]$ and is Lipschitz continuous. The
bandwidths $h_j$ converge to zero and $nh_j/(\log n) \rightarrow\infty
$ as $n \rightarrow\infty$.
\end{enumerate}
\end{assumptions*}

The assumption (A1) implies the concurvity condition (\ref
{concurvity}) since it implies that there exists a constant $c>0$ such
that for $\bff\in\cH(\bM)$,\vspace*{-1pt}
%
%
\begin{equation}\label{normbd}
\|\bff\|_\bM^2 \ge c \sum_{j=1}^d \int f_j(x_j)^2
p_j(x_j)\,\mathrm{d}x_j,\vspace*{-1pt}
\end{equation}
where $p_j$ denotes the marginal density function of $X_j$. The
inequality (\ref{normbd}) also tells us that the convergence of ${\hat
\mathbf{m}}$ in $\cH(\bM)$ implies the convergence of each component
$m_j$ in the usual $L_2$ norm.

\begin{theorem}\label{convthm}
Assume that \textup{(A1)}--\textup{(A5)} hold. Then, with probability tending to one,
there exists a solution $\{\hm_j\}_{j=1}^d$ of the backfitting
equation (\ref{backeqn1}) or (\ref{backeqn2}) that is unique.
Furthermore, there exist constants $0<\gamma<1$ and $0<C<\infty$ such
that, with probability tending to one,\vspace*{-1pt}
\[
\sum_{j=1}^d \int\bigl[\hm_j^{[r]}(x_j) - \hm_j(x_j) \bigr]^2
p_j(x_j) \,\mathrm{d}x_j \le C \gamma^{2r} \sum_{j=1}^d \int\bigl[\tm
_j(x_j)^2 + \hm_j^{[0]}(x_j)^2 \bigr] p_j(x_j) \,\mathrm{d}x_j.\vspace*{-3pt}
\]
\end{theorem}

\subsection{Asymptotic distribution of the backfitting estimators}\vspace*{-3pt}\label{sec3.2}

Next, we present the asymptotic distributions of $\hm_j$.
Define\vspace*{-1pt}
\begin{eqnarray*}
\tbm^A(\mathbf{x}) = \hbM(\mathbf{x})^{-1} n^{-1} \sum_{i=1}^n \bZ^i
[Y^i - m(\bX^i,\bZ^i) ]K_\bh(\mathbf{x}, \bX^i),\vspace*{-1pt}
\end{eqnarray*}
where $m(\bX,\bZ)$ is as given in (\ref{model}), and let $\tbm^B =
\tbm- \tbm^A$. As in the proof of Theorem~\ref{convthm}, we can
prove that, for $s=A$ or $B$, there exists a unique solution ${\hat
\mathbf{m}}^s \in\cH(\hbM)$ of the corresponding backfitting equation
(\ref{backeqn2}) where $\tbm$ is replaced by $\tbm^s$. By the
uniqueness of ${\hat\mathbf{m}}$, it follows that ${\hat\mathbf
{m}}= {\hat
\mathbf{m}}^A + {\hat\mathbf{m}}^B$.

Put ${\hat\mathbf{m}}^A = (\hm_1^A, \ldots, \hm_d^A)^\top\in\cH
(\hbM)$. In the proof of the following theorem, we will show that $\hm
_j^A$ are well approximated by $\tm_j^A \equiv(\hPi_j \tbm^A)_j$.
Note that
\[
(\hPi_j \tbm^A)_j(x_j) = {\hat{q}}_j(x_j)^{-1} n^{-1}\sum_{i=1}^n
Z_j^i [Y^i-m(\bX^i,\bZ^i) ]K_{h_j}(x_j,X_j^i).\vadjust{\goodbreak}
\]
Assume that the bandwidths $h_j$ are asymptotic to $c_j n^{-1/5}$ for
some constants $0<c_j<\infty$. By the standard techniques of kernel
smoothing, it can be proven that, for $\mathbf{x}$ in $(0,1)^d$, $
(\tm_1^A(x_j), \ldots, \tm_d^A(x_d) )^\top$, and thus ${\hat
\mathbf{m}}^A$, is asymptotically normal with mean zero and variance
$n^{-4/5}\operatorname{diag} (v_j(x_j) )$, where
\[
v_j(x_j) = \frac{E [Z_j^2 \sigma^2(\bX,\bZ) | X_j = x_j
]}{c_j p_j(x_j) [E (Z_j^2 | X_j = x_j) ]^2}\int
K(u)^2\,\mathrm{d}u
\]
and $\sigma^2(\bX,\bZ) = \operatorname{var} (Y | \bX, \bZ)$.
Here, it is worth noting that the vector $\tbm^A$, which belongs to
$L_2(\hbM)$, does not equal $ (\tm_1^A(x_j), \ldots, \tm
_d^A(x_d) )^\top\in\cH(\hbM)$.

The bias of the estimator ${\hat\mathbf{m}}$ comes from ${\hat
\mathbf{m}}^B$, which is the projection of $\tbm^B=(\tm_1^B, \ldots
, \tm
_d^B)^\top$ onto $\cH(\hbM)$. Define $\etam(\mathbf{x}) =
(c_1^2m_1''(x_1), \ldots, c_d^2 m_d''(x_d) )^\top$ and $\betam
_0 (\mathbf{x})$ by
\begin{eqnarray*}
\betam_0 (\mathbf{x}) &=& \Biggl[\sum_{k=1}^d c_k^2 m_k'(x_k) p(\mathbf
{x})^{-1}E(\bZ\bZ^\top|\bX=\mathbf{x})^{-1} \frac{\partial
}{\partial
x_k} \bigl(E(\bZ Z_k | \bX=\mathbf{x})p(\mathbf{x}) \bigr) + \frac{1}{2} \etam(\mathbf{x}) \Biggr]\\
&&  {}\times\int u^2 K(u)\,\mathrm{d}u.
\end{eqnarray*}
Note that $\tbm$ and $\betam_0$ do not belong to $\cH(\bM)$. In the
proof of the next theorem, we will show that $\betam_0 (\mathbf{x})$ is
the asymptotic bias of $\tbm(\mathbf{x})$ as an estimator of $\mathbf
{m}(\mathbf{x})$ and that the asymptotic bias of ${\hat\mathbf
{m}}(\mathbf{x})$
equals $\betam(\mathbf{x})$, where $\betam$ is the projection of
$\betam
_0$ onto $\cH(\bM)$:
\[
\betam\equiv\Pi(\betam_0 | \cH(\bM) ) =
\operatorname{argmin}\limits_{\bff\in\cH(\bM)} \int[\betam_0(\mathbf
{x}) - \bff
(\mathbf{x}) ]^\top\bM(\mathbf{x}) [\betam_0(\mathbf{x}) - \bff
(\mathbf{x}) ]\,\mathrm{d}\mathbf{x}.
\]
We write $\betam(\mathbf{x}) = (\beta_1(x_1), \ldots, \beta
_d(x_d) )^\top$.

The following theorem, which demonstrates the asymptotic joint
distribution of $\hm_j$, requires an additional condition on $m_j$.

\begin{enumerate}[(A6)]
\item[(A6)] $E(\bZ\bZ^\top\sigma^2(\bX,\bZ) | \bX=\mathbf{x})$ is
continuous on $[0,1]^d$.
\item[(A7)] The coefficient functions $m_j$ are twice continuously
differentiable on $[0,1]$, and $E(Z_j Z_k | \bX=\mathbf{x})$ is
continuously partially differentiable on $[0,1]^d$ for all $1 \le j, k
\le d$.
\end{enumerate}

\begin{theorem}\label{distthm}
Assume that \textup{(A1)--(A7)} hold and that the bandwidths $h_j$ are
asymptotic to $c_j n^{-1/5}$ for some constants $0<c_j<\infty$. Then,
for any $\mathbf{x}\in(0,1)^d$, $n^{2/5} [\hm_j(x_j) - m_j(x_j)]$
for $1
\le j \le d$ are jointly asymptotically normal with mean $ (\beta
_1(x_1), \ldots, \beta_d(x_d) )^\top$ and variance $\operatorname
{diag} (v_j(x_j) )$.
\end{theorem}

\section{The method with local polynomial fitting}\label{sec4}

The method we studied in the previous two sections is based on local
constant fitting, where we approximate $f_j(X_j^i)$ by $f_j(x_j)$ when
$X_j^i$ are near $x_j$, in the least-squares criterion $\sum_{i=1}^n
[Y^i - \sum_{j=1}^d f_j(X_j^i) Z_j^i ]^2$. The method may
be extended\vadjust{\goodbreak} to local polynomial fitting, where we approximate
$f_j(X_j^i)$ by $f_j(x_j) + (X_j^i-x_j)f_j^{(1)}(x_j) + \cdots+
(X_j^i-x_j)^\pi f_j^{(\pi)}(x_j)/\pi!$ for $X_j^i$ near $x_j$. Here
and below, $g^{(k)}$ denotes the $k$th derivative of a~function $g\dvtx
\mathbb{R} \rightarrow\mathbb{R}$. Define
\[
\bw_j(x_j,u_j) = \biggl(1, \biggl(\frac{u_j-x_j}{h_j} \biggr),
\ldots, \biggl(\frac{u_j-x_j}{h_j} \biggr)^\pi\biggr)^\top.
\]
We consider the following kernel-weighted least-squares criterion to
estimate $m_j$:
%
%
\begin{equation}\label{crt-pol}
L(\bff) = \int n^{-1} \sum_{i=1}^n \Biggl[Y^i - \sum_{j=1}^d Z_j^i
\bw_j(x_j, X_j^i)^\top\bff_j(x_j) \Biggr]^2 K_\bh(\mathbf{x}, \bX
^i)\,\mathrm{d}\mathbf{x},
\end{equation}
where $\bff^\top= (\bff_1^\top, \ldots, \bff_d^\top)$ and $\bff
_j(x_j) = (f_{j,0}(x_j), \ldots, f_{j,\pi}(x_j))^\top$ for functions
$f_{j,k}\dvtx  \mathbb{R} \rightarrow\mathbb{R}$. Let ${\hat\mathbf
{m}}$ be
the minimizer of $L(\bff)$. The proposed estimators of $m_j$ are then
$\hm_{j,0}$ in ${\hat\mathbf{m}}$, and~$\hm_{j,k}$ in~${\hat
\mathbf{m}}$
are estimators of $h^k m_j^{(k)}/k!$. We thus define the proposed
estimators of~$m_j^{(k)}$ by
\[
\hm_j^{(k)}(x_j) = k! h_j^{-k} \hm_{j,k}(x_j), \qquad  0 \le k \le\pi
,\ 1 \le j \le d.
\]

The minimization of $L(\bff)$ at (\ref{crt-pol}) is done over $\bff$
with $L(\bff) < \infty$. Define
\begin{eqnarray*}
\mathbf{v}(\bu,\mathbf{z};\mathbf{x})^\top&=& (\bw
_1(x_1,u_1)^\top z_1,
\ldots, \bw_d(x_d,u_d)^\top z_d ).
\end{eqnarray*}
The expression in the bracket at (\ref{crt-pol}) can then be written
as $Y^i- \mathbf{v}(\bX^i,\bZ^i;\mathbf{x})^\top\bff(\mathbf
{x})$. We now
redefine $\hbM$ used in the previous two sections as
%
%
\begin{equation}\label{def-Mhat}
\hbM(\mathbf{x}) = n^{-1}\sum_{i=1}^n \mathbf{v}(\bX^i,\bZ
^i;\mathbf{x})\mathbf{v}(\bX^i,\bZ^i;\mathbf{x})^\top K_\bh
(\mathbf{x}, \bX^i).
\end{equation}
$L(\bff) <\infty$ is then equivalent to $\int\bff(\mathbf{x})^\top
\hbM(\mathbf{x}) \bff(\mathbf{x})\,\mathrm{d}\mathbf{x}< \infty$
and minimizing
$L(\bff)$ is equivalent to minimizing $\int[\tbm(\mathbf{x})-
\bff(\mathbf{x}) ]^\top\hbM(\mathbf{x}) [\tbm(\mathbf{x})- \bff
(\mathbf{x}) ] \,\mathrm{d}\mathbf{x}$, where $\tbm(\mathbf{x})$
is redefined as
%
%
\begin{equation}\label{mtil-pol}
\tbm(\mathbf{x}) = \hbM(\mathbf{x})^{-1} n^{-1} \sum_{i=1}^n
\mathbf{v}(\bX
^i,\bZ^i;\mathbf{x}) Y^i K_\bh(\mathbf{x}, \bX^i).
\end{equation}
The function space that arises in this general problem is the class of
$(\pi+1)d$-vectors of functions $\bff= (f_{j,k})$ such that $\int
\bff(\mathbf{x})^\top\hbM(\mathbf{x}) \bff(\mathbf{x}) \,\mathrm
{d}\mathbf{x}< \infty$
and $f_{j,k}(\mathbf{x})=g_{j,k}(x_j)$ for some functions
$g_{j,k}\dvtx \mathbb
{R} \rightarrow\mathbb{R}$, $1 \le j \le d$ and $0 \le k \le\pi$.
We continue to denote the function space by $\cH(\hbM)$, and its norm
by $\|\cdot\|_\hbM$.
Thus,
%
%
\begin{equation}\label{def-est-pol}
{\hat\mathbf{m}}= \operatorname
{argmin}\limits_{\bff\in\cH(\hbM)} \|
\tbm- \bff\|_\hbM^2.
\end{equation}

By considering the G\^{a}teaux or Fr\'{e}chet derivatives of the
objective function $L(\bff)$ with respect to $\bff$, the solution
${\hat\mathbf{m}}$ of the minimization\vadjust{\goodbreak} problem (\ref{def-est-pol})
satisfies the following system of integral equations:
%
%
\begin{equation}\label{backeqn-pol}
\bzero= \int\hbM_{j}(\mathbf{x})^\top[\tbm(\mathbf{x}) - {\hat
\mathbf{m}}(\mathbf{x}) ] \,\mathrm{d}\mathbf{x}_{-j}, \qquad  1 \le j \le d,
\end{equation}
where $\bzero$ is the $(\pi+1)$-dimensional zero vector and $\hbM
_{j}$ are $(\pi+1)d \times(\pi+1)$ matrices defined by $\hbM=\hbM
^\top=(\hbM_{1}, \ldots, \hbM_{d})$.
We write ${\hat\mathbf{m}}^\top= ({\hat\mathbf{m}}_1^\top, \ldots
, {\hat
\mathbf{m}}_d^\top)$. Define
\begin{eqnarray*}
\tbm_
j(x_j) &=& \hPsim_j(x_j)^{-1}n^{-1}\sum_{i=1}^n \bw
_j(x_j,X_j^i) K_{h_j}(x_j,X_j^i) Z_j^i Y^i,\\
\hPsim_j(x_j) &=& n^{-1}\sum_{i=1}^n \bw_j(x_j,X_j^i)\bw
_j(x_j,X_j^i)^\top K_{h_j}(x_j,X_j^i) (Z_j^i)^2,\\
\hPsim_{jk}(x_j,x_k) &=& n^{-1}\sum_{i=1}^n \bw_j(x_j,X_j^i)\bw
_k(x_k,X_k^i)^\top K_{h_j}(x_j,X_j^i) K_{h_k}(x_k, X_k^i) Z_j^i Z_k^i
\end{eqnarray*}
for $k \neq j$. We then find that the system of $(\pi+1)$-dimensional
equations (\ref{backeqn-pol}) is equivalent to the following
backfitting equations which update the estimators of $m_j$ and their
derivatives up to the $\pi$th order:
%
%
\begin{equation}\label{backeqn1-pol}
{\hat\mathbf{m}}_j(x_j) = \tbm_j(x_j) - \sum_{k=1, \neq j}^d \int
\hPsim_j(x_j)^{-1}\hPsim_{jk}(x_j,x_k) {\hat\mathbf{m}}_k(x_k)\,
\mathrm{d}x_k, \qquad
1 \le j \le d.
\end{equation}

We want to emphasize again that the method with local polynomial
fitting does not require computation of the full-dimensional estimator
$\tbm(\mathbf{x})$ at (\ref{mtil-pol}). It only requires one- and
two-dimensional smoothing to compute $\tbm_j$, $\hPsim_j$ and $\hPsim
_{jk}$, and involves inversion of~$\hPsim_j$ only. Although $\hPsim
_j$ are $(\pi+1) \times(\pi+1)$ matrices, they are computed by means
of one-dimensional local smoothing so that they do not suffer from
sparsity of data in high dimensions. The marginal integration method,
in contrast, requires the computation of~$\tbm(\mathbf{x})$ and thus, in
practice, the marginal integration method may break down in the case
where $d$ is large.

\begin{BA*}
With a set of initial estimates ${\hat\mathbf{m}}_j^{[0]} = (\hm
_{j,0}, \ldots, \hm_{j,\pi} )^\top
$, we iterate for $r=1, 2, \ldots$ the following process: for $1 \le j
\le d$,
%
%
\begin{eqnarray}\label{backfit-pol}
{\hat\mathbf{m}}_j^{[r]}(x_j) &=& \tbm_j(x_j) - \sum_{k=1}^{j-1}
\int
\hPsim_j(x_j)^{-1}\hPsim_{jk}(x_j,x_k) {\hat\mathbf
{m}}_k^{[r]}(x_k)\,\mathrm{d}x_k \nonumber
\\[-8pt]
\\[-8pt]
&& {}- \sum_{k=j+1}^{d} \int\hPsim_j(x_j)^{-1}\hPsim
_{jk}(x_j,x_k) {\hat\mathbf{m}}_k^{[r-1]}(x_k)\,\mathrm{d}x_k.
\nonumber
\end{eqnarray}
\end{BA*}

In the following two theorems, we show that the backfitting algorithm
(\ref{backfit-pol}) converges to ${\hat\mathbf{m}}_j, 1 \le j \le d,$
at a geometric\vadjust{\goodbreak} rate and that ${\hat\mathbf{m}}_j, 1 \le j \le d,$ are
jointly asymptotically normal. We give the results for the case where
$\pi$, the order of local polynomial fitting, is odd. It is widely
accepted that fitting odd orders of local polynomial is better than
even orders. It also gives simpler formulas in the asymptotic expansion
and requires a weaker smoothness condition on $E(\bZ\bZ^\top| \bX
=\mathbf{x})$. In fact, instead of (A6) in Section~\ref{sec3}, we need the
following assumption:

\begin{enumerate}[(A7$'$)]
\item[(A7$'$)] The coefficient functions $m_j$ are $(\pi+1)$-times
continuously differentiable on $[0,1]$ and $E(Z_j Z_k | \bX=\mathbf{x})$
is continuous on $[0,1]^d$ for all $1 \le j, k \le d$.
\end{enumerate}

To state the first theorem, we need to introduce the limit of the
matrix $\hbM(\mathbf{x})$. Note that $\hbM(\mathbf{x})$ consists of
$(\pi
+1)\times(\pi+1)$ blocks\vspace*{-1pt}
\[
\hbM_{j,k}(\mathbf{x}) \equiv n^{-1}\sum_{i=1}^n \bw_j(x_j, X_j^i)
\bw
_k(x_k, X_k^i)^\top Z_j^i Z_k^i K_\bh(\mathbf{x}, \bX^i), \qquad  1 \le j,
k \le d.\vspace*{-1pt}
\]
For $j \neq k$, the matrices $\hbM_{j,k}(\mathbf{x})$ are
approximated by\vspace*{-1pt}
\[
E [\bw_j(x_j, X_j) \bw_k(x_k, X_k)^\top Z_j^i Z_k^i K_\bh(\mathbf
{x}, \bX) ] \simeq\bmu\bmu^\top E(Z_j Z_k | \bX=\mathbf{x})
p(\mathbf{x}),\vspace*{-1pt}
\]
where $\bmu= (\mu_\ell(K) )^\top$ and $\mu_\ell(K) =
\int u^\ell K(u)\,\mathrm{d}u$. On the other hand, for $j=k$,\vspace*{-1pt}
\[
\hbM_{j,j}(\mathbf{x}) \simeq\bN_1 E(Z_j^2 | \bX=\mathbf{x})
p(\mathbf{x}),\vspace*{-1pt}
\]
where $\bN_1$ is a $(\pi+1)\times(\pi+1)$ matrix defined by $\bN_1
= (\mu_{\ell+ \ell'}(K) )$. Here, we adopt the
convention that the indices of the matrix entries run from $(0,0)$ to
$(\pi,\pi)$. Thus, $\hbM(\mathbf{x})$ is approximated by\vspace*{-1pt}
%
%
\begin{equation}\label{def-M}
\bM(\mathbf{x}) \equiv p(\mathbf{x})\bigl [E(\bZ\bZ^\top| \bX=\mathbf{x})
\otimes(\bmu\bmu^\top) + \operatorname{diag} \bigl(E(Z_j^2 | \bX
=\mathbf{x})
\bigr) \otimes(\bN_1 - \bmu\bmu^\top) \bigr],\vspace*{-1pt}
\end{equation}
where $\otimes$ denotes the Kronecker product. The matrix $\bM
(\mathbf{x})$ is positive definite under the assumption (A1). To see
this, note
first that by (A1), the matrix $E(\bZ\bZ^\top| \bX=\mathbf{x})
\otimes
(\bmu\bmu^\top)$ is nonnegative definite. Also, \mbox{$E(Z_j^2 | \bX
=\mathbf{x})$} are bounded away from zero on $[0,1]^d$ for all $1 \le j
\le
d$. Furthermore, $\bN_1-\bmu\bmu^\top$ is the variance--covariance
matrix of $(1, U, \ldots, U^\pi)^\top$, where $U$ is a random
variable with density $K$. Since $K$ is supported on a uncountable set,
it follows that $\bN_1-\bmu\bmu^\top$ is positive definite. The
foregoing arguments show that the smallest eigenvalue of $\bM(\mathbf
{x})$ is bounded away from zero on $[0,1]^d$. Let $\cH(\bM)$ be
defined as $\cH(\hbM)$ with $\hbM$ being replaced by $\bM$ and
define its norm by $\|\bff\|_\bM^2 = \int\bff(\mathbf{x})^\top\bM
(\mathbf{x})\bff(\mathbf{x})\,\mathrm{d}\mathbf{x}$.\vspace*{-3pt}

\begin{theorem}\label{convthm-pol}
Assume that \textup{(A1)--(A5)} hold. Then, with probability tending to one,
there exists a solution $\{{\hat\mathbf{m}}_j\}_{j=1}^d$ of the
backfitting equation (\ref{backeqn1-pol}) that is unique. Furthermore,
there exist constants $0<\gamma<1$ and $0<C<\infty$ such that, with
probability tending to one,\vspace*{-1pt}
\begin{eqnarray*}
&&\sum_{j=1}^d \int\bigl|{\hat\mathbf{m}}_{j}^{[r]}(x_j) - {\hat\mathbf
{m}}_{j}(x_j) \bigr|^2 p_j(x_j) \,\mathrm{d}x_j \\[-2pt]
&& \quad  \le C \gamma^{2r} \sum_{j=1}^d \int\bigl[ |\tbm
_j(x_j) |^2 + \bigl|{\hat\mathbf{m}}_j^{[0]}(x_j) \bigr|^2 \bigr]
p_j(x_j) \,\mathrm{d}x_j.\vadjust{\goodbreak}
\end{eqnarray*}
\end{theorem}

In the next theorem, we give the asymptotic distribution of the
proposed estimators. We define $\mathbf{m}(\mathbf{x})^\top=
(\mathbf{m}_1(x_1)^\top, \ldots, \mathbf{m}_d(x_d)^\top)$, where
%
%
\begin{equation}\label{mult_m}
\mathbf{m}_j(x_j) = \bigl(m_j(x_j), h_j m_j^{(1)}(x_j)/1!, \ldots,
h_j^\pi m_j^{(\pi)}(x_j)/\pi! \bigr)^\top.
\end{equation}
For the bandwidths $h_j$, we assume that $h_j$ is asymptotic to $c_j
n^{-1/(2\pi+3)}$ for some constant $0<c_j<\infty$. Define $\gamm=
(\mu_{\pi+1}(K), \ldots, \mu_{\pi+1+\pi}(K) )^\top$
and a $(\pi+1)\times(\pi+1)$ matrix by $\bN_2 = (\mu_{\ell
+\ell'}(K^2) )$. For $1 \le j \le d$, define $\betam_j(x_j) =
c_j^{\pi+1}\bN_1^{-1} \gamm m_j^{(\pi+1)}(x_j)/ (\pi+1)!$ and
\[
\bV_j(x_j) = \frac{E [Z_j^2 \sigma^2(\bX,\bZ) | X_j =
x_j ]}{c_j p_j(x_j) [E (Z_j^2 | X_j = x_j) ]^2} \bN
_1^{-1} \bN_2 \bN_1^{-1}.
\]

\begin{theorem}\label{distthm-pol}
Assume that \textup{(A1)--(A6)} and \textup{(A7$'$)} hold, and that the bandwidths $h_j$
are asymptotic to $c_j n^{-1/(2\pi+3)}$ for some constants
$0<c_j<\infty$. Then, for any $\mathbf{x}\in(0,1)^d$, $n^{(\pi
+1)/(2\pi
+3)}\times  [{\hat\mathbf{m}}_j(x_j) - \mathbf{m}_j(x_j) ]$, $1 \le j
\le d$, are asymptotically independent and
\[
n^{(\pi+1)/(2\pi+3)} [{\hat\mathbf{m}}_j(x_j) - \mathbf{m}_j(x_j)
] \overset{d}{\rightarrow} N (\betam_j(x_j),
\bV_j(x_j) ), \qquad  1 \le j \le d.
\]
\end{theorem}

Theorem~\ref{distthm-pol} not only gives the asymptotic distributions
of the estimators of the coefficient functions $m_j$, but also those of
their derivatives. Recall the definition of $\mathbf{m}_j$ at (\ref
{mult_m}) and that ${\hat\mathbf{m}}_j(x_j) = (\hm_j(x_j), h_j \hm
_j^{(1)}(x_j)/1!, \ldots, h_j^\pi\hm_j^{(\pi)}(x_j)/\pi!
)^\top$. Thus, the theorem implies that $n^{(\pi+1-k)/(2\pi+3)}
[\hm_j^{(k)}(x_j) - m_j^{(k)}(x_j) ]$ is asymptotically
normal with mean $k! c_j^{\pi+1-k} (\bN_1^{-1} \gamm)_k\times  m_j^{(\pi
+1)}(x_j) /(\pi+1)!$ and variance
\[
(k!)^2 (\bN_1^{-1} \bN_2 \bN_1^{-1})_{kk} \frac{E [Z_j^2
\sigma^2(\bX,\bZ) | X_j = x_j ]}{c_j^{2k+1}
p_j(x_j) [E (Z_j^2 | X_j = x_j) ]^2},
\]
where, for a vector $\mathbf a$ and a matrix $\mathbf B$, ${\mathbf
a}_k$ denotes
the $k$th entry of $\mathbf a$ and ${\mathbf B}_{kk}$ denotes the $k$th
diagonal entry of $\mathbf B$. In the case of local linear fitting
($\pi
=1$), we have
\[
(\bN_1^{-1}\bN_2 \bN_1^{-1})_{00} = \int K^2(u)\,\mathrm{d}u, \qquad  (\bN
_1^{-1}\gamm)_0 = \int u^2 K(u)\,\mathrm{d}u.
\]

Another implication of Theorem~\ref{distthm-pol} is that the
estimators $\hm_j^{(k)}(x_j)$ for $0 \le k \le\pi$ have the oracle
properties. Suppose that we know all other coefficient functions except
$m_j$. In this case, we would estimate $m_j$ and its derivatives up to
order $\pi$ by minimizing
\[
n^{-1} \sum_{i=1}^n \Biggl[Y^i - \sum_{k=1, \neq j}^d m_k(X_k^i)Z_k^i
- Z_j^i \bw_j(x_j, X_j^i)^\top\bff_j(x_j) \Biggr]^2 K_{h_j}(x_j, X_j^i)
\]
over $\bff_j$. It can be shown that the resulting estimators of
$m_j^{(k)}(x_j)$ for $0 \le k \le\pi$ have\vspace*{-2pt} the same asymptotic
distributions as $\hm_j^{(k)}(x_j)$ for $0 \le k \le\pi$.

The asymptotically optimal choices of the bandwidths $h_j$ may be
derived from Theorem~\ref{distthm-pol}. Let $c_j^{\pi+1} b_{j}(x_j)$
and $c_j^{-1} \tau_{j}(x_j)$ denote the asymptotic mean and the
asymptotic variance of $n^{(\pi+1)/(2\pi+3)} [\hm_j(x_j) -
m_j(x_j) ]$, respectively. That is,
\begin{eqnarray*}
b_j(x_j)&=&(\bN_1^{-1} \gamm)_0 m_j^{(\pi+1)}(x_j)/ (\pi+1)!,\\
\tau_j(x_j)&=& \frac{E [Z_j^2 \sigma^2(\bX,\bZ) | X_j =
x_j ]}{p_j(x_j) [E (Z_j^2 | X_j = x_j) ]^2} (\bN
_1^{-1} \bN_2 \bN_1^{-1})_{00}.
\end{eqnarray*}
The optimal choice of $c_j$ which minimizes the asymptotic mean
integrated squared error is then given by
%
%
\begin{equation}\label{opt-band}
c_j^{\mathrm{opt}} = \biggl[\frac{\int\tau_j(x_j) p_j(x_j)\,\mathrm
{d}x_j}{2(\pi
+1)\int b_j(x_j)^2 p_j(x_j)\,\mathrm{d}x_j} \biggr]^{1/(2\pi+3)}.
\end{equation}
This formula for the optimal bandwidth involves unknown quantities. We
may get a~rule-of-thumb bandwidth selector by fitting polynomial
regression models, as in Yang \textit{et al.}~\cite{Yang2006}, to estimate the unknown
quantities in the formula for $c_j^{\mathrm{opt}}$; see Section~\ref{sec6},
where we
employ this approach to analyze Boston Housing Data. Alternatively, we
may adopt the approach of Mammen and Park \cite{Mammen2005} to obtain more
sophisticated bandwidth selectors.

\section{Numerical properties}\label{sec5}

We investigated the finite-sample properties of the proposed estimators
in comparison with the marginal integration method studied in Yang \textit{et al.}
\cite{Yang2006}. We considered the case where local linear
smoothing ($\pi
=1$) is employed. The simulation study was done in two settings, one in
a low-dimensional case ($d=3$) and the other in a high-dimensional case
($d=10$).

In the first case, we generated the data $(\bX^i, \bZ^i, Y^i)$ from
the model:
$Y = m_1 (X_1)Z_1 + m_2(X_2)Z_2 + m_3(X_3)Z_3 + \sigma(\bX,\bZ
)\varepsilon$, where $Z_1 \equiv1$ and
%
%
\begin{equation}\label{sigma-sim}
\sigma(\mathbf{x},\mathbf{z}) = \frac{1}{2} + \frac
{z_2^2+z_3^2}{1+z_2^2+z_3^2}\exp\biggl(-2+\frac{x_1+x_2}{2} \biggr).
\end{equation}
The vector $\bX=(X_1, X_2, X_3)$ was generated from the uniform
distribution over the unit cube $(0,1)^3$, and the covariate vector
$(Z_2, Z_3)$ was generated from the bivariate normal with mean $(0,0)$
and covariance matrix
$\left({1 \atop 0.5 }\enskip{ 0.5 \atop 1}\right)$.
The vectors $\bX$ and $\bZ$ were independent, and the error
term $\varepsilon$ was generated from the standard normal
distribution, independently of $(\bX,\bZ)$. This model was also
considered in Yang \textit{et al.} \cite{Yang2006}. We took $m_1(x)=1+\mathrm{e}^{2x-1}$,
$m_2(x)=\cos(2\uppi x)$ and $m_3(x)=x^2$.

In the second case, where $d=10$, we took the same variance function
$\sigma^2(\mathbf{x},\mathbf{z})$ as in~(\ref{sigma-sim}), for the
sake of
simplicity. Thus, $\sigma^2(\mathbf{x},\mathbf{z})$ did not depend
on $(x_j,
z_j)$ for $4 \le j \le10$. The extra covariates $X_j$ for $4 \le j \le
10$ were generated from the uniform distribution over $(0,1)^7$
independently of $(X_1, X_2, X_3)$, and $Z_j$ for $4 \le j \le10$ were
generated from the multivariate normal distribution with mean $\bzero$
and covariance matrix $\bI$, the identity matrix, independently of
$(Z_2, Z_3)$ and of $\bX$. We chose $m_j(x) = x^2$ for $4 \le j
\le10$.\vadjust{\goodbreak}

We used the Epanechnikov kernel $K(u) = (3/4)(1- u^2)I[-1,1](u)$ and
the optimal bandwidths $h_j^{\mathrm{sbf}}=c_j^{\mathrm
{opt}}n^{-1/5}$, where
$c_j^{\mathrm{opt}}$ are given at (\ref{opt-band}). This was for the
proposed estimator. For the marginal integration method, the estimator
$\hat{m}_j^{\mathrm{mi}}$ of the $j$th coefficient function $m_j$
that we
investigated was
\[
\hat{m}_j^{\mathrm{mi}}(x_j) = n^{-1}\sum_{i=1}^n \hat\theta_j(X_1^i,
\ldots, X_{j-1}^i, x_j, X_{j+1}^i, \ldots, X_j^n),
\]
where $\hat\theta_j(\mathbf{x})$ was the $[(j-1)(\pi+1)+1]$st entry of
$\tbm(\mathbf{x})$ defined at (\ref{mtil-pol}), but $K_{h_k}$ (for $k
\neq j$) in the definition of $\tbm(\mathbf{x})$ was replaced by
$L_{b_k}$. Note that, for the marginal integration method, in the
estimation of the $j$th coefficient function, we may choose another
kernel $L$ and need to use other bandwidths $b_k$, different from
$h_j$, for the directions of $x_k (k \neq j)$ not of interest. We took
$L=K$ and $b_k = c(\log n)^{-1}h_j^{\mathrm{mi}}$ for all directions
$k \neq
j$, where $h_j^{\mathrm{mi}}$ is the optimal bandwidth for the marginal
integration method, obtained similarly as the one for the proposed
method at (\ref{opt-band}), and $c$ was a constant multiplier for
which we tried four values, $c=1, 3, 5, 10$.

We used $\tbm_j$ defined in Section~\ref{sec4} as the initial estimates ${\hat
\mathbf{m}}_j^{[0]}$ for the proposed method. The backfitting algorithm
converged very fast. We took
\[
\sqrt{\sum_{j=1}^d \int\bigl[\hm_j^{[r-1]}(x_j) - \hm
_j^{[r]}(x_j) \bigr]^2\,\mathrm{d}x_j} \le 10^{-11}
\]
as a criterion for the convergence. With this criterion, the
backfitting algorithm converged within 11 iterations in all cases. The
average number of iterations was 6.5 from the 500 replications.
In a preliminary numerical study with the marginal integration method,
we found that inverting the matrix $\hbM(\mathbf{x})$ often caused
numerical instability of the estimates, even for the low-dimensional
case where $d=3$. This reflects the curse of dimensionality that the
marginal integration suffers from. Thus, we actually computed a
`ridged' version of $\hat{m}_j^{\mathrm{mi}}$ by adding $n^{-2}$ to the
diagonal entries of the matrix $\hbM(\mathbf{x})$. The same modification
was also made in the numerical study of Yang \textit{et al.} \cite{Yang2006}.

%
\begin{table*}
\tabcolsep=0pt
\caption{The mean integrated squared errors (MISE), the integrated
squared biases (ISB) and the integrated variances (IV) of the marginal
integration estimators (MI) and the proposed estimators (SBF) when
$d=3$ (the constant $c$ for MI is the multiplier $c$ in the formula
$b_k = c (\log n)^{-1} h_j^{\mathrm{mi}}$, where $b_k$ is the secondary
bandwidth applied to the direction of $x_k$, $k \neq j$, in the
estimation of $m_j$)}
\label{tab1}\begin{tabular*}{\textwidth}{@{\extracolsep{\fill}}llllllll@{}}
\hline
\multirow{2}{27pt}{Sample size}& \multirow{2}{40pt}{Coefficient function}& & \multicolumn{4}{c}{MI} &SBF\\[-5pt]
&  & & \multicolumn{4}{c}{\hrulefill} &\\
 &  & & $c=1$ & $c=3$ & $c=5$ & $c=10$ & \\
\hline
$n=100$ & $m_1$ & MISE & 0.1190 & 0.1140 & 0.1158 & 0.1151 & 0.1496 \\
& & ISB & 0.0174 & 0.0150 & 0.0147 & 0.0145 & 0.0019 \\
& & IV & 0.1016 & 0.0990 & 0.1011 & 0.1006 & 0.1476 \\
& $m_2$ & MISE & 0.6354 & 0.5738 & 0.5795 & 0.5826 & 0.3613 \\
& & ISB & 0.4089 & 0.3502 & 0.3465 & 0.3461 & 0.0484 \\
& & IV & 0.2265 & 0.2236 & 0.2330 & 0.2364 & 0.3129 \\
& $m_3$ & MISE & 0.1873 & 0.2218 & 0.2255 & 0.2259 & 0.2512 \\
& & ISB & 0.0057 & 0.0056 & 0.0056 & 0.0056 & 0.0017 \\
& & IV & 0.1816 & 0.2163 & 0.2200 & 0.2203 & 0.2495 \\
[6pt]
$n=400$ & $m_1$ & MISE & 0.0347 & 0.0332 & 0.0365 & 0.0363 & 0.0415 \\
& & ISB & 0.0092 & 0.0087 & 0.0087 & 0.0086 & 0.0005 \\
& & IV & 0.0255 & 0.0245 & 0.0279 & 0.0277 & 0.0410 \\
& $m_2$ & MISE & 0.2648 & 0.2815 & 0.2872 & 0.2894 & 0.1244 \\
& & ISB & 0.2126 & 0.2227 & 0.2248 & 0.2257 & 0.0199 \\
& & IV & 0.0521 & 0.0588 & 0.0624 & 0.0637 & 0.1045 \\
& $m_3$ & MISE & 0.0478 & 0.0576 & 0.0610 & 0.0620 & 0.0810 \\
& & ISB & 0.0050 & 0.0046 & 0.0046 & 0.0047 & 0.0008 \\
& & IV & 0.0428 & 0.0529 & 0.0564 & 0.0573 & 0.0802 \\
\hline
\end{tabular*}\vspace*{-0.5pt}
\end{table*}

Table~\ref{tab1} shows the results for the case $d=3$, based on 500 data sets
with sizes $n=100$ and $400$. The table provides the mean integrated
squared errors (MISE) of the estimators of each coefficient function
$m_j$, defined by
\begin{eqnarray*}
\operatorname{MISE}_j (\bar{m}_j) &=& \int E [\bar{m}_j(x_j) -
m_j(x_j) ]^2 \,\mathrm{d}x_j\\[2pt]
&=& \int[E \bar{m}_j(x_j) - m_j(x_j) ]^2 \,\mathrm{d}x_j + \int
[\bar{m}_j(x_j) - E \bar{m}_j(x_j) ]^2 \,\mathrm{d}x_j\\[2pt]
&\overset{\mathrm{let}}{=}& \operatorname{ISB}_j (\bar{m}_j) +
\operatorname{IV}_j (\bar{m}_j)
\end{eqnarray*}
for an estimator $\bar{m}_j$. It also gives the integrated squared
bias (ISB) and the integrated variance (IV). The results suggest that
the proposed method gives better performance in terms\vadjust{\goodbreak} of $\operatorname
{MISE}_{\mathrm{tot}} = \sum_{j=1}^{3} \operatorname{MISE}_j$. When
$n=100$, the sum
of $\operatorname{MISE}_j$ of $\hm_j$ equals $0.7621$, while those of the
marginal integration method are $0.9417, 0.9096, 0.9208,
0.9236$ for $c=1, 3, 5, 10$, respectively. In the case where
$n=400$, $\operatorname{MISE}_{\mathrm{tot}}=0.2469$ for the proposed
method, while
it equals $0.3473, 0.3723, 0.3847, 0.3877$ for the marginal
integration method.

According to Table~\ref{tab1}, the performance of the marginal integration
method appears not to be sensitive to the choice of the secondary
bandwidth $b_k$. However, this is true only when we use the optimal
bandwidth $h_j^{\mathrm{mi}}$. In fact, we found that the performance
depended crucially on the choice $b_k$ when other choices of $h_j$ were
used. As an example, we report in Table~\ref{tab2} the results when one uses
$h_j=h_j^{\mathrm{mi}}/3$ instead of $h_j=h_j^{\mathrm{mi}}$. In the
latter case,
the sum of $\operatorname{MISE}_j$ ranges from $0.8001$ to $2.7453 $ when
$n=100$, and from $0.2291$ to $2.1080$ when $n=400$, for those four
values of $c$. One interesting thing to note is that the ISB of the
marginal integration increases drastically as $c$ decreases. The main
lesson here is that the choice of the secondary bandwidths $b_k$ for
the marginal integration method is as important as the choice of $h_j$.

The finite-sample results in Table~\ref{tab1} show some discrepancy with the
asymptotics for the functions $m_1$ and $m_3$. Asymptotically, if the
optimal bandwidth is used, then the IV is four times as large as the
ISB. In general, finite-sample properties do not always match with
asymptotics. One possible reason for the discrepancy in this particular
setting is that the coefficient functions $m_1$ and $m_3$ are far
simpler than the complexity brought by the noise level, so the proposed
method easily catches the structure with less bias. This seems not to
be the case with the marginal integration, however. For the marginal
integration, the secondary bandwidths $b_k$ interact with the primary
bandwidth $h_j$ for the bias and variance performance, as discussed in
the previous paragraph.

%
\begin{table*}
\tabcolsep=0pt
\caption{The mean integrated squared errors (MISE), the integrated
squared biases (ISB) and the integrated variances (IV) of MI when
$h_j=h_j^{\mathrm{mi}}/3$ was used (the constant $c$ is the multiplier $c$
in the formula $b_k = c (\log n)^{-1} h_j^{\mathrm{mi}}$, where $b_k$
is the
secondary bandwidth applied to the direction of $x_k$, $k \neq j$, in
the estimation of $m_j$)}
\label{tab2}\begin{tabular*}{\textwidth}{@{\extracolsep{\fill}}lllllll@{}}
\hline
\multirow{2}{27pt}{Sample size}&\multirow{2}{60pt}{Coefficient function} & & \multicolumn{4}{c@{}}{MI} \\[-5pt]
 & & & \multicolumn{4}{c@{}}{\hrulefill}\\
& & & $c=1$ & $c=3$ & $c=5$ & $c=10$ \\
\hline
$n=100$ & $m_1$ & MISE & 1.5109 & 0.2664 & 0.1822 & 0.1737 \\
& & ISB & 1.4327 & 0.0096 & 0.0011 & 0.0012 \\
& & IV & 0.0782 & 0.2568 & 0.1812 & 0.1725 \\
& $m_2$ & MISE & 0.6611 & 0.4578 & 0.3459 & 0.3576 \\
& & ISB & 0.3095 & 0.0340 & 0.0338 & 0.0313 \\
& & IV & 0.3516 & 0.4238 & 0.3121 & 0.3263 \\
& $m_3$ & MISE & 0.5733 & 0.2743 & 0.2720 & 0.3012 \\
& & ISB & 0.0217 & 0.0013 & 0.0014 & 0.0015 \\
& & IV & 0.5516 & 0.2730 & 0.2706 & 0.2997 \\
[6pt]
$n=400$ & $m_1$ & MISE & 1.4539 & 0.0891 & 0.0465 & 0.0465 \\
& & ISB & 1.4177 & 0.0032 & 0.0004 & 0.0003 \\
& & IV & 0.0362 & 0.0859 & 0.0461 & 0.0462 \\
& $m_2$ & MISE & 0.3554 & 0.1596 & 0.1109 & 0.1129 \\
& & ISB & 0.2359 & 0.0139 & 0.0154 & 0.0160 \\
& & IV & 0.1195 & 0.1457 & 0.0955 & 0.0969 \\
& $m_3$ & MISE & 0.2987 & 0.0702 & 0.0717 & 0.0856 \\
& & ISB & 0.0188 & 0.0005 & 0.0006 & 0.0007 \\
& & IV & 0.2799 & 0.0697 & 0.0711 & 0.0849 \\
\hline
\end{tabular*}
\vspace*{-3pt}
\end{table*}

%
\begin{table*}
\tabcolsep=0pt
\caption{The mean integrated squared errors (MISE), the integrated
squared biases (ISB) and the integrated variances (IV) of the marginal
integration estimators (MI) and the proposed estimators (SBF) when $d=10$}
\label{tab3}\begin{tabular*}{\textwidth}{@{\extracolsep{\fill}}llllllll@{}}
\hline
\multirow{2}{27pt}{Sample size}& \multirow{2}{40pt}{Coefficient function}& \multicolumn{3}{c}{MI}&\multicolumn{3}{c@{}}{SBF}\\[-5pt]
 &  &
\multicolumn{3}{c}{\hrulefill}&\multicolumn{3}{c@{}}{\hrulefill}\\
 &  & MISE & ISB & IV & MISE & ISB & IV \\
\hline
$n=100$ & $m_1$ & 0.2533 & 0.1242 & 0.1291 & 0.1904 & 0.0046 & 0.1858 \\
& $m_2$ & 0.7284 & 0.4353 & 0.2931 & 0.4357 & 0.0605 & 0.3752 \\
& $m_3$ & 0.2622 & 0.0059 & 0.2563 & 0.3042 & 0.0024 & 0.3018 \\
& $m_4$ & 0.1303 & 0.0054 & 0.1249 & 0.1404 & 0.0022 & 0.1382 \\
& $m_5$ & 0.1351 & 0.0060 & 0.1291 & 0.1489 & 0.0011 & 0.1478 \\
& $m_6$ & 0.1336 & 0.0055 & 0.1281 & 0.1509 & 0.0019 & 0.1490 \\
& $m_7$ & 0.1345 & 0.0054 & 0.1291 & 0.1677 & 0.0019 & 0.1658 \\
& $m_8$ & 0.1228 & 0.0053 & 0.1175 & 0.1482 & 0.0019 & 0.1463 \\
& $m_9$ & 0.1428 & 0.0071 & 0.1357 & 0.1707 & 0.0009 & 0.1698 \\
& $m_{10}$ & 0.1270 & 0.0059 & 0.1211 & 0.1528 & 0.0014 & 0.1514 \\
[6pt]
$n=400$ & $m_1$ & 0.0505 & 0.0115 & 0.0390 & 0.0457 & 0.0008 & 0.0449 \\
& $m_2$ & 0.2999 & 0.2223 & 0.0776 & 0.1264 & 0.0196 & 0.1068 \\
& $m_3$ & 0.0642 & 0.0054 & 0.0588 & 0.0893 & 0.0004 & 0.0889 \\
& $m_4$ & 0.0324 & 0.0048 & 0.0276 & 0.0379 & 0.0005 & 0.0374 \\
& $m_5$ & 0.0358 & 0.0054 & 0.0304 & 0.0355 & 0.0010 & 0.0345 \\
& $m_6$ & 0.0369 & 0.0040 & 0.0329 & 0.0331 & 0.0004 & 0.0327 \\
& $m_7$ & 0.0300 & 0.0044 & 0.0256 & 0.0370 & 0.0009 & 0.0361 \\
& $m_8$ & 0.0319 & 0.0043 & 0.0276 & 0.0368 & 0.0006 & 0.0362 \\
& $m_9$ & 0.0321 & 0.0052 & 0.0269 & 0.0364 & 0.0009 & 0.0355 \\
& $m_{10}$ & 0.0303 & 0.0046 & 0.0257 & 0.0355 & 0.0006 & 0.0349 \\
\hline
\end{tabular*}
\vspace*{-3pt}
\end{table*}

Table~\ref{tab3} shows the results for the case $d=10$. Here, for the marginal
integration, we report only the results when $c=5$ which gave the best
performance. In fact, the marginal integration got worse very quickly
as $c$ decreased from $c=5$. For example, we found the total MISE,
$\sum_{j=1}^{10} \operatorname{MISE}_j$, was $3.4996$ when $c=3$ and was
$6.1834$ when $c=1$, in the case where $n=400$. Note that the value
equals $0.6440$ when $c=5$ and $n=400,$ as reported in Table~\ref{tab3}. For the
proposed method, it equals $0.5136$.

\section{Analysis of Boston Housing Data}\label{sec6}

The data consist of fourteen variables, among which one is response and
the other thirteen are predictors. There are $506$ observations from
506 tracts in the Boston area; see Harrison and Rubinfeld \cite
{Harrison1978} for
details about the data set. The data set has been analyzed by Fan and
Huang \cite{Fan2005} and Wang and Yang \cite{Wang2009}, among
others.\vadjust{\goodbreak}
The former fitted
the data using a partially linear functional coefficient model where
all coefficient functions in the nonparametric part are functions of a
single variable. The latter considered an additive regression model.
Here, we apply the varying coefficient model (\ref{model}) to fit the
data using the proposed method. We take the variable MEDV (median value
of owner-occupied homes in \$1000's) as the response variable $Y$. We
consider five variables as covariates $X_j$ or $Z_j$. They are CRIM
(per capita crime rate by town), RM (average number of rooms per
dwelling), TAX (full-value property tax rate per \$10\,000),
PTRATIO (pupil--teacher ratio by town) and LSTAT (percentage of lower
income status of the population). As in Wang and Yang \cite{Wang2009},
we take
logarithmic transformation for TAX and LSTAT to remove sparse areas in
the domains of these variables.

We want to find a varying coefficient model that fits the data set
well. Since LSTAT can be a good explanatory variable that determines
the overall level of the housing price, we consider models of the form
%
%
\begin{equation}\label{modelgp}
\operatorname{MEDV} = m_1(\log(\operatorname{LSTAT})) + m_2(X_2)Z_2
+ m_3(X_3)Z_3 +\mathrm{(noise)}.\vadjust{\goodbreak}
\end{equation}
A general question is which variables should be the model covariates
$Z_j$ and which should take the role of $X_j$. This may be obvious for
some data sets, but it is not so clear for the Boston Housing Data.
Thus, we fitted all possible models and chose the one that best fitted
the data. In general, we do not suggest employing the
all-possible-models approach since it can get out of control quickly as
the number of variables increases, and it induces a certain
arbitrariness in the choice. For the Boston Housing Data, there are
only twelve varying coefficient models of the form (\ref{modelgp}),
listed in Table~\ref{tab4}, and all models are interpretable. If the number of
variables is large, then we suggest first choosing a~set of model
covariates $Z_j$ among all covariates by fitting parametric linear
models and using a variable selection technique, and then picking one
as $X_j$ for each $Z_j$ from the remaining variables based on a
criterion such as RSPE (which is defined later).

%
\begin{table*}
\tabcolsep=0pt
\caption{Relative squared prediction errors obtained from fitting 12
varying coefficient models with the Boston Housing Data}
\label{tab4}\begin{tabular*}{\textwidth}{@{\extracolsep{\fill}}llllll@{}}
\hline
\multirow{2}{45pt}{Model no.}& \multicolumn{4}{c}{Covariates} &
\multirow{2}{70pt}{Relative squared prediction error} \\[-5pt]
 & \multicolumn{4}{c}{\hrulefill} & \\
 & $X_2$ & $Z_2$ & $X_3$ & $Z_3$ & \\
\hline
\hphantom{1}1 & CRIM & RM & TAX & PTRATIO & 0.3514 \\
\hphantom{1}2 & CRIM & RM & PTRATIO & TAX & N/A \\
\hphantom{1}3 & CRIM & TAX & RM & PTRATIO & 0.2700 \\
\hphantom{1}4 & CRIM & TAX & PTRATIO & RM & 0.2688 \\
\hphantom{1}5 & CRIM & PTRATIO & RM & TAX & 0.4390 \\
\hphantom{1}6 & CRIM & PTRATIO & TAX & RM & 0.4757 \\
\hphantom{1}7 & RM & CRIM & TAX & PTRATIO & 0.3010 \\
\hphantom{1}8 & RM & CRIM & PTRATIO & TAX & \textit{0.2412} \\
\hphantom{1}9 & RM & TAX & PTRATIO & CRIM & N/A \\
10 & RM & PTRATIO & TAX & CRIM & N/A \\
11 & TAX & CRIM & PTRATIO & RM & N/A \\
12 & TAX & RM & PTRATIO & CRIM & N/A \\
\hline
\end{tabular*}
\end{table*}

We employed local linear smoothing in implementing the proposed method
and used the Epanechnikov kernel. For the bandwidths $h_j$, we chose to
use a rule-of-thumb method that we describe below. Note that the
unknowns in the expression of the optimal bandwidth at (\ref
{opt-band}) are $A_j = \int m_j''(x_j)^2 p_j(x_j)\,\mathrm{d}x_j$,
$B_j(x_j) =
E [Z_j^2 \sigma^2(\bX,\bZ) | X_j =x_j]$ and $C_j(x_j) = E
(Z_j^2 | X_j = x_j)$.
The second derivative of $m_j$ in $A_j$ can be estimated by fitting a
cubic polynomial regression model. This gives $\hat A_j = n^{-1}\sum
_{i=1}^n (2 \hat\alpha_{j,2} + 6 \hat\alpha_{j,3}X_j^i)^2$,
where $\hat\alpha_{j,k}$ are the least-squares estimators that minimize
\[
\sum_{i=1}^n \Biggl[Y^i - \sum_{j=1}^d (\alpha_{j,0} + \alpha
_{j,1}X_j^i + \alpha_{j,2}X_j^{i 2} + \alpha_{j,3}X_j^{i 3}
)Z_j^i \Biggr]^2.
\]
Here, we take $Z_1^i \equiv1$. The conditional means, $B_j$ and $C_j$,
can be estimated by fitting linear regression models. Since $E[Z_j^2
\sigma^2(\bX,\bZ) | X_j = x_j] = E[Z_j^2 (Y-m(\bX,\bZ))^2
| X_j=x_j]$, the conditional mean $B_j$ is estimated by $\hat
B_j(x_j) = \hat\beta_{j,0} + \hat\beta_{j,1} x_j$, where $\hat
\beta_{j,0}$ and $\hat\beta_{j,1}$ minimize
\[
\sum_{i=1}^n \Biggl[Z_j^{i 2} \Biggl(Y^i - \sum_{k=1}^d (\hat
\alpha_{k,0} + \hat\alpha_{k,1}X_k^i + \hat\alpha_{k,2}X_k^{i 2}
+ \hat\alpha_{k,3}X_k^{i 3} )Z_k^i \Biggr)^2 - \beta_{j,0} -
\beta_{j,1} X_j^i \Biggr]^2.
\]
Similarly, $C_j$ for $j=2, 3$ are estimated by $\hat C_j(x_j) = \hat
\gamma_{j,0} + \hat\gamma_{j,1} x_j$, where $\hat\gamma_{j,0}$ and
$\hat\gamma_{j,1}$ minimize $\sum_{i=1}^n (Z_j^{i 2} - \gamma
_{j,0} - \gamma_{j,1} X_j^i)^2$. Note that $C_1 \equiv1$.

We split the data set into two parts, one for estimation of the models
and the other for assessment of the estimated models. We selected 100
tracts for the model assessment out of 506 distributed in 92 towns.
This was done in a manner that would lead to more selections in a town
with a larger number of tracts. We fitted the twelve varying
coefficient models using the data for the remaining 406 tracts and made
out-of-sample predictions with the data for the selected 100 tracts. We
calculated their relative squared prediction errors,
\[
\operatorname{RSPE} = \frac{\sum_{i=1}^{100} [\operatorname{MEDV}^i
- \hat
m_1(\log(\operatorname{LSTAT}^i)) - \hat m_2(X_2^i)Z_2^i - \hat
m_3(X_3^i)Z_3^i ]^2}{\sum_{i=1}^{100} [\operatorname{MEDV}^i -
\overline{\operatorname{MEDV}} ]^2},
\]
where $\hat m_j$ for $j=1,2,3$ were constructed by using the data for
the 406 remaining tracts.

Table~\ref{tab4} reports the results. In the table, we do not provide the values
of RSPE for the models numbered 2, 9, 10, 11 and 12. In the preliminary
fitting of these models taking~$X_j$ and $Z_j$ as specified, we found
that they produced extremely large residuals for some of the
observations that corresponded to $\operatorname{PTRATIO} = 20.2$ or
$\operatorname{TAX}
= 666$. This resulted in a negative value of $\hat B_j (x_j)$ for a
certain range of $x_j$ and, as a consequence, produced a~negative
estimate of $\int\tau_j(x_j)p_j(x_j)\,\mathrm{d}x_j$ in the
bandwidth formula
(\ref{opt-band}). Since these five models do not explain MEDV well as
a function of the covariates and would give a~large value of RSPE when
fitted, we excluded them from further analysis.

\begin{figure}

\includegraphics{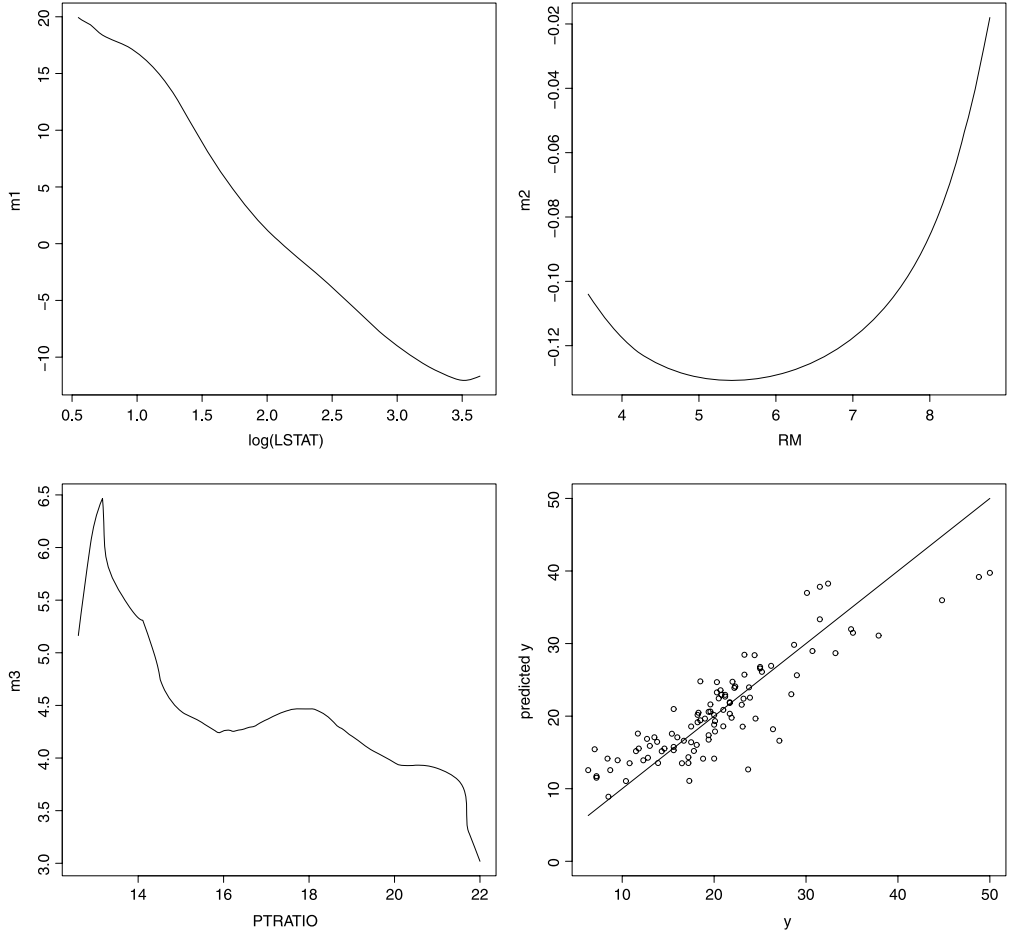}

\caption{For the final model (\protect\ref{modelfit}), the
upper-left, upper-right and lower-left panels depict the estimated
coefficient functions $\hat m_1$, $\hat m_2$ and $\hat m_3$,
respectively, and the lower-right panel exhibits plots of the observed
values $Y^i$ versus their predicted values $\hat Y^i$.}
\label{fig:brainslices}
\label{fig1}\end{figure}

According to the table, the model with the smallest RSPE is
%
%
\begin{equation}\label{modelfit}
\operatorname{MEDV} = m_1(\log(\operatorname{LSTAT})) +
m_2(\operatorname{RM})\operatorname{CRIM} {}+{}
m_3(\operatorname{PTRATIO})\log(\operatorname{TAX}) + \mathrm{(noise)}.
\end{equation}
Figure~\ref{fig1} depicts the estimated coefficient functions $\hat m_1, \hat
m_2$ and $\hat m_3$. It also plots the actual values of MEDV and their
predicted values according to the estimated model from~(\ref
{modelfit}). The prediction was made for those 100 tracts that were not
used in estimating the model. The estimated curve $\hat m_1$ indicates
that a high percentage of lower income status decreases the prices of
homes. The estimated curve $\hat m_2$ suggests that for towns with
higher or lower average numbers of rooms per dwelling, the crime rate
is less influential on the prices of homes. Finally, from the estimated
curve $\hat m_3$, we see that
if the pupil--teacher ratio gets higher, then the prices of homes
increase less rapidly as the property tax rate increases. The curve
$\hm_3$ looks somewhat rigid. The reason for this is that the variable
PTRATIO does not really take values on a continuous scale since it is
the pupil--teacher ratio by town, so that all tracts in a town have the
same value of PTRATIO. Furthermore, some towns share the same value
with others. For example, the 132 tracts (out of 506) associated with
the 15 towns in the city of Boston have the same value, $20.2$.

%

\begin{appendix}\label{appm}
\setcounter{equation}{0}
\section*{Appendix: Technical details}

\subsection{\texorpdfstring{Proof of Theorem~\protect\ref{convthm}}{Proof of Theorem 1}}\label{sec6.1}
We prove that there exists a constant $0<\gamma<1$ such that $\|\hQ\|
< \gamma$ with probability tending to one. Let $\cH_j(\bM)$ be
defined as $\cH_j(\hbM)$ with $\hbM$ being replaced by $\bM$. Let
$p_j$ and $p_{jk}$ denote the marginal densities of $X_j$ and $(X_j,
X_k)$, respectively. Define
%
%
\begin{eqnarray}\label{condz}
q_j(x_j) &=& E(Z_j^2 |X_j = x_j) p_j(x_j),\\\label{condz1}
q_{jk}(x_j, x_k) &=& E(Z_j Z_k | X_j=x_j, X_k=x_k) p_{jk}(x_j,x_k),
\qquad
k \neq j.\vadjust{\goodbreak}
\end{eqnarray}
For $\bff_j \in\cH_j(\bM)$,
%
%
\begin{eqnarray}\label{normcal}
\|\bff_j\|_\bM^2 = \int\bff_j(\mathbf{x})^\top\bM(\mathbf{x})
\bff
_j(\mathbf{x})\,\mathrm{d}\mathbf{x}
= \int f_j(x_j)^2 q_j(x_j)\,\mathrm{d}x_j.
\end{eqnarray}
The equality (\ref{normcal}) follows from the identity
\[
\int E(Z_j^2 | \bX=\mathbf{x}) p(\mathbf{x}) \,\mathrm{d}\mathbf
{x}_{-j} = E(Z_j^2
|X_j=x_j) p_j(x_j).
\]
From (\ref{normcal}) and H\"{o}lder's inequality, it follows that, for
$\bff\in\cH(\bM)$,
\begin{eqnarray*}
&&\|(\hQ_j - Q_j)\bff\|_\bM\\[2pt]
&& \quad  = \biggl[\int\biggl(\sum_{k=1,\neq j} \int\biggl[\frac
{{\hat{q}}_{jk}(x_j,x_k)}{{\hat{q}}_j(x_j)} - \frac
{q_{jk}(x_j,x_k)}{q_j(x_j)} \biggr]f_k(x_k)\,\mathrm{d}x_k \biggr)^2 q_j(x_j)
\,\mathrm{d}x_j \biggr]^{1/2}\\[2pt]
&& \quad  \le\sum_{k=1, \neq j} \biggl[\int\biggl(\frac{{\hat
{q}}_{jk}(x_j,x_k)}{{\hat{q}}_j(x_j)q_k(x_k)} - \frac
{q_{jk}(x_j,x_k)}{q_j(x_j)q_k(x_k)} \biggr)^2 q_j(x_j)q_k(x_k)\,\mathrm{d}x_j
\,\mathrm{d}x_k \biggr]^{1/2}\\[2pt]
&&\hphantom{\le\sum_{k=1, \neq j}} \quad {} \times\biggl[\int f_k(x_k)^2 q_k(x_k)\,\mathrm{d}x_k
\biggr]^{1/2}\\[2pt]
&& \quad  \le \mathrm{o}_p(1) \sum_{k=1, \neq j} \|\bff_k\|_\bM.
\end{eqnarray*}
Since $\|Q_j\| =1$, this proves that $\|\hQ_j\| \le C_1$ with
probability tending to one for some constant $0<C_1<\infty$. Define $Q
= Q_d \cdots Q_1$. Then,
\begin{eqnarray*}
\|\hQ- Q\| = \Biggl\|\sum_{k=0}^{d-1}Q_d \cdots Q_{d-k+1}(\hQ_{d-k} -
Q_{d-k}) \hQ_{d-k-1}\cdots\hQ_1 \Biggr\| = \mathrm{o}_p(1),
\end{eqnarray*}
where we interpret both $Q_{d+1}$ and $\hQ_0$ as the zero operator.
From (A1), (A3) and (\ref{normcal}), the projection operators $\Pi_j\dvtx
\cH_k(\bM) \rightarrow\cH_j(\bM)$ for all $1 \le j \neq k \le d$
are Hilbert--Schmidt. By applying parts B, C and D of Proposition A.4.2
of Bickel, Klaassen, Ritov and Wellner~\cite{Bickel1993}, we find that
$\|Q\| <1$.
This shows that there exists a constant $0<\gamma<1$ such that $\|\hQ
\| <\gamma$ with probability tending to one.

To complete the proof of Theorem~\ref{convthm}, it follows from (\ref
{backfit2}) that with probability tending to one,
\[
\bigl\|{\hat\mathbf{m}}^{[r]} - {\hat\mathbf{m}}\bigr\|_\bM= \Biggl\|\sum
_{s=r}^\infty\hQ^s \hbr+ \hQ^r {\hat\mathbf{m}}^{[0]} \Biggr\|_\bM
\le \gamma^r \biggl(\|\hbr\|_\bM\frac{1}{1-\gamma} + \bigl\|{\hat
\mathbf{m}}^{[0]}\bigr\|_\bM\biggr).
\]
By (\ref{backfit3}) and the fact that $\|\hQ_j\| \le C_1$ with
probability tending to one, there exists a constant $0<C_2<\infty$
such that with probability tending to one,
\[
\|\hbr\|_\bM\le C_2 \sum_{j=1}^d\biggl [\int\tm_j(x_j)^2 q_j(x_j)
\,\mathrm{d}x_j \biggr]^{1/2}.
\]
This completes the proof of Theorem~\ref{convthm}.

\subsection{\texorpdfstring{Proof of Theorem~\protect\ref{distthm}}{Proof of Theorem 2}}\label{sec6.2}

We will prove that for each $\mathbf{x}\in(0,1)^d$,\vspace*{-1pt}
%
%
\begin{eqnarray}\label{th2:pf1}
\hm_j^A(x_j) &=& \tm_j^A(x_j) + \mathrm{o}_p(n^{-2/5}) \qquad  \mbox{for } 1
\le j \le d,\\\label{th2:pf3}
{\hat\mathbf{m}}^B(\mathbf{x}) &=& \mathbf{m}(\mathbf{x}) + \betam
(\mathbf{x})
n^{-2/5} + \mathrm{o}_p(n^{-2/5}).
\end{eqnarray}

\begin{pf*}{Proof of (\ref{th2:pf1})} Note that ${\hat\mathbf{m}}^A = \sum
_{s=0}^\infty\hQ^s \hbr^A$, where\vspace*{-1pt}
%
%
\begin{equation}\label{th2:pf1.5}
\hbr^A = (I-\hQ)\tbm^A = \tbm_d^A + \hQ_d \tbm_{d-1}^A + \cdots+
\hQ_d \cdots\hQ_2 \tbm_1^A
\end{equation}
and $\tbm_j^A(\mathbf{x}) = (0, \ldots, 0, \tm_j^A(x_j), 0, \ldots,
0)^\top$. From formulas (\ref{form1})--(\ref{proj}), it follows that
\[
\hQ_d \cdots\hQ_{j+1} \tbm_j^A (\mathbf{x}) = (0, \ldots, 0,
\tm_j^A(x_j), \tg_{j+1}(x_{j+1}), \ldots, \tg_d(x_d) )^\top
, \qquad  2 \le j \le d,\
\]
for some random functions $\tg_k\dvtx  \mathbb{R} \rightarrow\mathbb
{R}$, $j+1 \le k \le d$, where the first $j-1$ entries of the vector on
the right-hand side of the equation are zero. This implies that
%
%
\begin{equation}\label{th2:pf4}
\hbr^A(\mathbf{x}) = (\tm_1^A(x_1), \hg_2(x_2), \ldots, \hg
_d(x_d) )^\top,
\end{equation}
where $\hg_k$ for $2 \le k \le d$ are random functions from $\mathbb
{R}$ to $\mathbb{R}$. If we prove
that
%
%
\begin{equation}\label{th2:pf5}
\sup_{\mathbf{x}\in[0,1]^d} \Biggl|\sum_{s=1}^\infty\hQ^s \hbr^A
(\mathbf{x}) \Biggr| = \mathrm{o}_p(n^{-2/5}),
\end{equation}
then (\ref{th2:pf4}) implies (\ref{th2:pf1}) for the case $j=1$. By
exchanging the entries of $\tbm^A$, we can see that (\ref{th2:pf1})
also holds for $j \ge2$.

To prove (\ref{th2:pf5}), it suffices to show that
%
%
\begin{eqnarray}\label{th2:pf6}
\sup_{\mathbf{x}\in[0,1]^d} |\hQ\hbr^A (\mathbf{x}) |&=&
\mathrm{o}_p(n^{-2/5}),\\\label{th2:pf7}
\Biggl\|\sum_{s=1}^\infty\hQ^s \hbr^A \Biggr\|_\bM&=&
\mathrm{o}_p(n^{-2/5}).
\end{eqnarray}
To see this, note that from (\ref{form1}) and (\ref{form2}), we have,
for $\bff= (f_1, \ldots, f_d)^\top\in\cH(\hbM)$,
%
%
\begin{equation}\label{th2:pf8}
\hQ_j \bff(\mathbf{x}) = (f_1(x_1), \ldots, f_{j-1}(x_{j-1}),
f_j^*(x_j), f_{j+1}(x_{j+1}), \ldots, f_d(x_d) )^\top,\vadjust{\goodbreak}
\end{equation}
where $f_j^*(x_j) = - \sum_{k=1, \neq j}^d \int f_k(x_k) \frac{{\hat
{q}}_{jk}(x_j,x_k)}{{\hat{q}}_j(x_j)}\,\mathrm{d}x_k$.
Thus, there exists a constant \mbox{$0<C<\infty$} such that with probability
tending to one,
\[
\sup_{\mathbf{x}\in[0,1]^d} \Biggl|\sum_{s=2}^\infty\hQ^s \hbr^A
(\mathbf{x}) \Biggr| = \sup_{\mathbf{x}\in[0,1]^d} \Biggl| \hQ\sum
_{s=1}^\infty
\hQ^s \hbr^A (\mathbf{x}) \Biggr|\le C \Biggl\|\sum_{s=1}^\infty\hQ^s
\hbr^A \Biggr\|_\bM.
\]
We prove (\ref{th2:pf6}) and (\ref{th2:pf7}). From standard kernel
theory, we can prove that for all $k \neq j$,
%
%
\begin{equation}\label{th2:pf9}
\sup_{x_k \in[0,1]} \biggl| \int\tm_j^A (x_j)\frac{{\hat
{q}}_{jk}(x_j,x_k)}{{\hat{q}}_k(x_k)}\,\mathrm{d}x_j \biggr| = \mathrm{o}_p(n^{-2/5}).
\end{equation}
The approximation (\ref{th2:pf9}), together with the expressions at
(\ref{th2:pf1.5}) and (\ref{th2:pf8}), gives (\ref{th2:pf6}). Since
$\|\hQ\| < \gamma$ with probability tending to one for some $0<\gamma
<1$, we have
\[
\Biggl\|\sum_{s=1}^\infty\hQ^s \hbr^A \Biggr\|_\bM\le\sum
_{s=2}^\infty\gamma^s \|\hQ\hbr^A \|_\bM= \mathrm{o}_p(n^{-2/5}).
\]
This completes the proof of (\ref{th2:pf1}).
\end{pf*}

\begin{pf*}{Proof of (\ref{th2:pf3})} Let $\bl
_1(\mathbf{x}, \bu) = ((u_1-x_1)m_1'(x_1), \ldots,
(u_d-x_d)m_d'(x_d) )^\top$ and $\bl_2(\mathbf{x},\bu)=
((u_1-x_1)^2m_1''(x_1)/2, \ldots, (u_d-x_d)^2m_d''(x_d)/2
)^\top$. To get an idea of which terms in an expansion of $\tbm
^B(\mathbf{x})$ lead to the main terms in the expansion (\ref{th2:pf3}),
we note from an expansion of $m(\bX^i)$ that $\tbm^B(\mathbf{x})$ is
approximated by
%
%
\begin{eqnarray}\label{th2:pf11}
&& \mathbf{m}(\mathbf{x}) + \hbM(\mathbf{x})^{-1}n^{-1}\sum
_{i=1}^n \bZ^i \bZ
^{i\top} \bl_1(\mathbf{x}, \bX^i) K_\bh(\mathbf{x}, \bX^i)\nonumber
\\[-8pt]
\\[-8pt]
&& \quad {} + \hbM(\mathbf{x})^{-1}n^{-1}\sum_{i=1}^n \bZ^i
\bZ^{i\top} \bl_2(\mathbf{x}, \bX^i) K_\bh(\mathbf{x}, \bX
^i).
\nonumber
\end{eqnarray}
Define $\tbm^{B,1}(\mathbf{x}) = \hbM(\mathbf{x})^{-1}\int\bM
(\mathbf{x})
\bl_1(\mathbf{x},\bu) K_\bh(\mathbf{x},\bu) \,\mathrm{d}\bu$.
The second term of
(\ref{th2:pf11}) is then approximated by $\tbm^{B,1}(\mathbf{x}) +
\bM
(\mathbf{x})^{-1}\sum_{k=1}^d [\partial\bM_k(\mathbf{x})/\partial
x_k ] h_k^2 m_k'(x_k) \int u^2 K(u)\,\mathrm{d}u$.
Also, the third term is approximated by $ (h_1^2 m_1''(x_1)/2,
\ldots, h_d^2 m_d''(x_d)/2 )^\top\int u^2 K(u) \,\mathrm{d}u$. Define
\begin{eqnarray*}
\tbm^{B,2}(\mathbf{x}) &=& \Biggl[\bM(\mathbf{x})^{-1}\sum_{k=1}^d \frac
{\partial}{\partial x_k} \bM_k(\mathbf{x}) h_k^2 m_k'(x_k) + \frac
{1}{2} (h_1^2 m_1''(x_1), \ldots, h_d^2 m_d''(x_d) )^\top
\Biggr]\\
&& {}\times\int u^2 K(u) \,\mathrm{d}u
\end{eqnarray*}
and let $\tbm^{B,3}(\mathbf{x}) = \tbm^B(\mathbf{x}) - \mathbf
{m}(\mathbf{x}) -
\tbm^{B,1}(\mathbf{x})- \tbm^{B,2}(\mathbf{x})$.

For $\ell=1,2,3$, define ${\hat\mathbf{m}}^{B,\ell}$ to be the solution
of the backfitting equation at (\ref{backeqn2}) with $\tbm$ being
replaced by $\tbm^{B,\ell}$. By arguing as in the proof of (\ref
{th2:pf1}), we can deduce that $\hm_j^{B,3}(x_j) = \mathrm{o}_p(n^{-2/5})$ for
all $x_j \in(0,1)$. The projection of $\tbm^{B,2}$ onto $\cH(\hbM)$
is well approximated by the projection onto $\cH(\bM)$ with a
remainder $\deltam$ such that $\deltam(\mathbf{x}) = \mathrm{o}_p(n^{-2/5})$ for
all $\mathbf{x}\in(0,1)^d$. This proves that ${\hat\mathbf
{m}}^{B,2}(\mathbf{x}) = \betam(\mathbf{x})n^{-2/5} + \mathrm{o}_p(n^{-2/5})$
for all $\mathbf{x}\in(0,1)^d$.

It thus remains to prove that ${\hat\mathbf{m}}^{B,1}(\mathbf{x}) =
\mathrm{o}_p(n^{-2/5})$ for all $\mathbf{x}\in(0,1)^d$. For this bound, we will
show that $\hm_j^{B,1}(x_j) = \mu_j(x_j) + \mathrm{o}_p(n^{-2/5})$, uniformly
for all $x_j \in[0,1]$, $1 \le j \le d$, where
$\mu_j(x_j) = a_j(x_j) / \int K_{h_j}(x_j,u_j) \,\mathrm{d}u_j$ and $a_j(x_j)
= m_j'(x_j) \int(u_j-x_j) K_{h_j}(x_j,u_j) \,\mathrm{d}u_j $. For a
proof of
this claim, it suffices to show that
%
%
\begin{equation} \label{addeq}
\int\hbM_j(\mathbf{x})^\top[\tbm^{B,1}(\mathbf{x}) - {\bolds{\mu
}}(\mathbf{x}) ]\,\mathrm{d}\mathbf{x}_{-j} = \mathrm{o}_p(n^{-2/5}),
\end{equation}
uniformly for all $x_j \in[0,1]$, $1 \le j \le d$. Here, ${\bolds{\mu
}}(\mathbf{x})= (\mu_1(x_1),\ldots,\mu_d(x_d))^{\top}$.

We prove (\ref{addeq}). Note that, uniformly for $x_j \in[0,1]$,
\begin{eqnarray*}
&&\int\hbM_j(\mathbf{x})^\top{\bolds{\mu}}(\mathbf{x}) \,\mathrm
{d}\mathbf{x}_{-j}\\
&& \quad  = \biggl[\int q_j(u_j)K_{h_j}(x_j,u_j)\,\mathrm{d}u_j \biggr] \mu
_j(x_j) \\
&& \qquad {} + \sum_{k=1, \neq j} \int\mu_k(x_k)\biggl [\int
q_{jk}(u_j, u_k) K_{h_j}(x_j,u_j)K_{h_k}(x_k,u_k)\,\mathrm{d}u_j\, \mathrm{d}u_k
\biggr]\,\mathrm{d}x_k \\
&& \qquad {} + \mathrm{o}_p(n^{-2/5})\\
&& \quad  = q_j(x_j) a_j(x_j) + \sum_{k=1, \neq j} \int a_k(x_k)
q_{jk}(x_j, x_k) \,\mathrm{d}x_k \int K_{h_j}(x_j,u_j) \,\mathrm
{d}u_j + \mathrm{o}_p(n^{-2/5}).
\end{eqnarray*}
Claim (\ref{addeq}) now follows from the fact that
\begin{eqnarray*}
&&\int\hbM_j(\mathbf{x})^\top\tbm^{B,1}(\mathbf{x}) \,\mathrm
{d}\mathbf{x}_{-j}\\
&& \quad  = q_j(x_j) a_j(x_j) + \sum_{k=1, \neq j} \int a_k(x_k)
q_{jk}(x_j, x_k) \,\mathrm{d}x_k \int K_{h_j}(x_j,u_j) \,\mathrm
{d}u_j + \mathrm{o}_p(n^{-2/5})
\end{eqnarray*}
uniformly for $x_j \in[0,1]$.
\end{pf*}

\subsection{\texorpdfstring{Proofs of Theorems~\protect\ref{convthm-pol} and \protect\ref{distthm-pol}}
{Proofs of Theorems 3 and 4}}\label{sec6.3}

Recall the definitions of $\hbM$ and $\bM$ at (\ref{def-Mhat}) and
(\ref{def-M}), respectively, in the case of local polynomial fitting.
Let $\cH_j(\hbM)$ denote the space of $(\pi+1)d$-vectors of
functions $\bff= (f_{j,k})$ in $L_2(\hbM)$ such that $f_{j,\ell
}(\mathbf{x})=g_{j,\ell}(x_j)$, $0 \le\ell\le\pi,$ for some functions
$g_{j,\ell}\dvtx  \mathbb{R} \rightarrow\mathbb{R}$ and $\bff_k \equiv
(f_{k,0}, \ldots, f_{k,\pi})^\top= \bzero$ for $k \neq j$. As in
the case of local constant fitting, we can write $\cH(\hbM) = \cH
_1(\hbM) + \cdots+ \cH_d (\hbM)$. Define $\cH_j(\bM)$ likewise.
The vectors of functions that take the roles of $q_j$ and $q_{jk}$,
respectively, are
\begin{eqnarray*}
\Psim_j(x_j) &=& \bN_1 E(Z_j^2 |X_j = x_j) p_j(x_j),\\[3pt]
\Psim_{jk}(x_j, x_k) &=& \bmu\bmu^\top E(Z_j Z_k | X_j=x_j, X_k=x_k)
p_{jk}(x_j,x_k), \qquad  k \neq j.
\end{eqnarray*}
We then have projection formulas analogous to (\ref{form1})--(\ref
{proj}). For example, for $\bff\in L_2(\hbM)$ and $\bg\in L_2(\bM
)$, we obtain
\begin{eqnarray*}
(\hPi_j \bff)_j &=& \hPsim_j(x_j)^{-1} \int\hbM_j(\mathbf
{x})^\top
\bff(\mathbf{x}) \,\mathrm{d}\mathbf{x}_{-j},\\
(\Pi_j \bg)_j &=& \Psim_j(x_j)^{-1} \int\bM_j(\mathbf{x})^\top
\bg
(\mathbf{x}) \,\mathrm{d}\mathbf{x}_{-j}
\end{eqnarray*}
and $(\hPi_j \bff)_k = \bzero= (\Pi_j \bg)_k$ for $k \neq j$,
where $(\hPi_j \bff)_k$ and $(\Pi_j \bg)_k$ denote the $k$th $(\pi
+1)$-vector of the projection of $\bff$ onto $\cH_j(\hbM)$ and of
$\bg$ onto $\cH_j(\bM)$, respectively.
We can proceed as in the proof of Theorem~\ref{convthm} to prove
Theorem~\ref{convthm-pol}.

We prove Theorem~\ref{distthm-pol}. Decompose $\tbm$ at (\ref
{mtil-pol}) as $\tbm^A + \tbm^B$, where
\[
\tbm^A(\mathbf{x}) = \hbM(\mathbf{x})^{-1} n^{-1} \sum_{i=1}^n
\mathbf{v}(\bX
^i,\bZ^i;\mathbf{x}) [Y^i-m(\bX^i, \bZ^i) ] K_\bh(\mathbf{x},
\bX^i).
\]
Define ${\hat\mathbf{m}}^A$ and ${\hat\mathbf{m}}^B$ from $\tbm^A$ and
$\tbm^B$, respectively, to be the solutions of the backfitting
equation (\ref{backeqn1-pol}). It follows that $(\hPi_j \tbm^A)_j
(x_j)=\tbm_j^A(x_j)$, where
\[
\tbm_j^A(x_j) = \hPsim_j(x_j)^{-1} n^{-1}\sum_{i=1}^n \bw
_j(x_j,X_j^i) K_{h_j}(x_j, X_j^i) Z_j^i [Y^i-m(\bX^i,\bZ
^i) ].
\]
As in the proof of Theorem~\ref{distthm}, we can prove that ${\hat
\mathbf{m}}_j^A(x_j) = \tbm_j^A(x_j)+\mathrm{o}_p(n^{-(\pi+1)/(2\pi+3)})$ for
all $\mathbf{x}\in(0,1)^d$. The stochastic term $\tbm_j^A(x_j)$ has mean
zero and is asymptotically normal. Since $\hPsim_j(x_j) = \Psim
_j(x_j) + \mathrm{o}_p(1)$ and
\begin{eqnarray*}
&&n^{-1}h_j\sum_{i=1}^n \operatorname{var} [\bw
_j(x_j,X_j^i)K_{h_j}(x_j,X_j^i)Z_j^i Y^i | \bX^i, \bZ^i
]\\
&& \quad  = \bN_2 E [Z_j^2 \sigma^2(\bX,\bZ) | X_j=x_j
]p_j(x_j) + \mathrm{o}_p(1),
\end{eqnarray*}
we find that the asymptotic variance of $\tbm_j^A(x_j)$ equals
\begin{eqnarray*}
n^{-1}h_j^{-1} (\bN_1^{-1}\bN_2 \bN_1^{-1} )\frac{E
[Z_j^2 \sigma^2(\bX,\bZ) | X_j = x_j
]}{p_j(x_j) [E (Z_j^2 | X_j = x_j) ]^2}.
\end{eqnarray*}

Next, we approximate ${\hat\mathbf{m}}^B(\mathbf{x})$. Define
\begin{eqnarray*}
\tbm^{B,1}(\mathbf{x}) &=& \frac{1}{(\pi+1)!} \bM(\mathbf{x})^{-1}
n^{-1}\sum_{i=1}^n \mathbf{v}(\bX^i, \bZ^i;\mathbf{x})\\[-2pt]
&&{} \times \Biggl[\sum_{j=1}^d Z_j^i \biggl(\frac
{X_j^i-x_j}{h} \biggr)^{\pi+1}m_j^{(\pi+1)}(x_j)h_j^{\pi+1}
\Biggr]K_\bh(\mathbf{x}, \bX^i)
\end{eqnarray*}
and $\tbm^{B,2}(\mathbf{x}) = \tbm^B(\mathbf{x}) - \mathbf
{m}(\mathbf{x}) - \tbm
^{B,1}(\mathbf{x})$. As in the proof of Theorem~\ref{distthm}, we can show that
${\hat
\mathbf{m}}_j^{B,2}(x_j)= \mathrm{o}_p(n^{-(\pi+1)/(2\pi+3)})$ for all $x_j
\in
(0,1)$. We compute ${\hat\mathbf{m}}^{B,1}(\mathbf{x})$. We can
prove that,
for all $x_j \in(0,1)$,
%
%
\begin{eqnarray}\label{int-eqn1}
&&\int\hbM_j(\mathbf{x})^\top\tbm^{B,1}(\mathbf{x})\,\mathrm
{d}\mathbf{x}_{-j}
\nonumber\\[-2pt]
&& \quad  = \frac{1}{(\pi+1)!} \biggl[\bmu\mu_{\pi+1} \sum_{k=1,
\neq j} \int q_{jk}(x_j,x_k)h_k^{\pi+1} m_k^{(\pi+1)}(x_k)\,\mathrm
{d}x_k \\[-2pt]
&& \quad \hphantom{= \frac{1}{(\pi+1)!} \biggl[} {}+ h_j^{\pi+1}\gamm q_j(x_j) m_j^{(\pi+1)}(x_j)
\biggr] + \mathrm{o}_p\bigl(n^{-(\pi+1)/(2\pi+3)}\bigr),\nonumber
\end{eqnarray}
where $q_j$ and $q_{jk}$ are as defined at (\ref{condz}) and (\ref
{condz1}), respectively, and $\mu_{\pi+1}=\mu_{\pi+1}(K)$.
We also have
%
%
\begin{eqnarray}\label{int-eqn2}
\int\hbM_j(\mathbf{x})^\top{\hat\mathbf{m}}^{B,1}(\mathbf{x})\,
\mathrm{d}\mathbf{x}_{-j} &=& \bmu\bmu^\top\sum_{k=1, \neq j} \int
q_{jk}(x_j,x_k){\hat\mathbf{m}}_k^{B,1}(x_k)\,\mathrm{d}x_k \nonumber
\\[-9pt]
\\[-9pt]
&&{} + \bN_1 q_j(x_j){\hat\mathbf{m}}_j^{B,1}(x_j) +
\mathrm{o}_p\bigl(n^{-(\pi+1)/(2\pi+3)}\bigr)
\nonumber
\end{eqnarray}
for all $x_j \in(0,1)$. Now, we observe that $\bmu^\top\bN
_1^{-1}=(1,0, \ldots, 0)$ since $\bmu$ is the first column of $\bN
_1$. Thus,
\[
\bmu\bmu^\top\bN_1^{-1} \gamm= \bmu(1, 0, \ldots, 0)\gamm= \bmu
\mu_{\pi+1}.
\]
Comparing the two systems of equations (\ref{int-eqn1}) and (\ref
{int-eqn2}), and by the uniqueness of~${\hat\mathbf{m}}^{B,1}$, we
conclude that
\[
{\hat\mathbf{m}}_j^{B,1}(x_j) = (\bN_1^{-1} \gamm) h_j^{\pi+1}
m_j^{(\pi+1)}(x_j) /(\pi+1)! + \mathrm{o}_p\bigl(n^{-(\pi+1)/(2\pi+3)}\bigr)
\]
for all $x_j \in(0,1)$, $1 \le j \le d$. This completes the proof of
Theorem~\ref{distthm-pol}.
\end{appendix}

\section*{Acknowledgements}
Y.K. Lee was supported by National Research Foundation of Korea (NRF)
Grant NRF-2010-616-C00008. E. Mammen was supported by the DFG-NRF
Cooperative Program. B.U. Park was supported by NRF Grant
No. 20100017437, funded by the Korea government (MEST).\vadjust{\goodbreak}

\printhistory

\end{document}